\documentclass[aos,MSNbibl,nameyear,dvips]{arximspdf}
\usepackage{graphicx}

\doi{10.1214/11-AOS909}
\volume{39}
\issue{6}
\pubyear{2011}
\firstpage{2795}
\lastpage{2819}

\makeatletter

\newcommand{\R}{{\mathbb R}}

\newcommand{\cov}{\operatorname{cov}}
\newcommand{\corr}{\operatorname{corr}}
\renewcommand{\hat}{\widehat}

\newcommand{\ep}{\varepsilon}
\newcommand{\la}{\lambda}

\newtheorem{theorem}{Theorem}
\newproclaim{conjecture}{Conjecture}
\newproclaim{assumption}{Assumption}
\newtheorem{lemma}{Lemma}

\newcommand{\erfc}{\operatorname{erfc}}
\newcommand{\var}{\operatorname{var}}

\makeatother

\begin{document}
\begin{frontmatter}
\vspace*{12pt}
\dochead{2010 Rietz lecture}
\title{When does the screening effect hold?}
\runtitle{The screening effect}

\begin{aug}
\author[A]{\fnms{Michael L.} \snm{Stein}\corref{}\thanksref{t1}\ead[label=e1]{stein@galton.uchicago.edu}}
\runauthor{M. L. Stein}
\affiliation{University of Chicago}
\address[A]{Department of Statistics\\
University of Chicago\\
Chicago, Illinois 60637\\
USA\\
\printead{e1}} 
\end{aug}

\thankstext{t1}{Supported by US Department of Energy Grant
DE-SC0002557.}

\received{\smonth{11} \syear{2010}}
\revised{\smonth{4} \syear{2011}}

%
\begin{abstract}
When using optimal linear prediction to interpolate point observations
of a mean square continuous stationary spatial process, one often finds
that the interpolant mostly depends on those observations located
nearest to the predictand. This phenomenon is called the screening
effect. However, there are situations in which a screening effect does
not hold in a reasonable asymptotic sense, and theoretical support for
the screening effect is limited to some rather specialized settings for
the observation locations. This paper explores conditions on the
observation locations and the process model under which an asymptotic
screening effect holds. A series of examples shows the difficulty in
formulating a general result, especially for processes with different
degrees of smoothness in different directions, which can naturally
occur for spatial-temporal processes. These examples lead to a~general
conjecture and two special cases of this conjecture are proven. The key
condition on the process is that its spectral density should change
slowly at high frequencies. Models not satisfying this condition of
slow high-frequency change should be used with caution.
\end{abstract}

%

\end{frontmatter}

\section{Introduction}\label{sec1}
The screening effect is the geostatistical term for the phenomenon of
nearby observations tending to reduce the influence of more distant
observations when using kriging (optimal linear prediction) for spatial
interpolation [\citet{JouHui78}, \citet{ChiDel99}]. This
phenomenon is often invoked as a justification for ignoring more
distant observations when using kriging [\citet{MemMou07},
\citet{Eme09}]. Only in some very limited special cases is the
effect exact in the sense\vadjust{\goodbreak} that the more distant observations make no
contribution to the kriging predictor, so it is natural to use
asymptotics as a way to study the screening effect.

Let us set some notation. Write $x\cdot y$ for the inner product of
commensurate vectors $x$ and $y$. Suppose $Z$ is a mean square
continuous, stationary, mean 0 Gaussian process on $\R^d$ with
autocovariance function $K(x) = E\{Z(x)Z(0)\}$ and spectral density
$f$, so that $K(x) = \int_{\R^d} e^{i\omega\cdot x} f(\omega) \,d\omega$.
When the mean is assumed known to be 0, kriging is often called simple
kriging. Throughout this work, we assume that the problem of interest
is to predict $Z(0)$. For $S\subset\R^d$, write $Z(S)$ for the vector
of observations (in some order) of $Z$ on~$S$, and define $e(S)$ to be
the error of the best linear predictor, or BLP, of $Z(0)$ based on
$Z(S)$. Let $N_\ep$ and $F_\ep$ be two classes of sets indexed by the
parameter $\ep>0$, with $N_\ep$ representing observations near 0 and
$F_\ep$ more distant observations. We will say that $N_\ep$
asymptotically screens out the effect of $F_\ep$ if
%
%
\begin{equation}\label{main}
\lim_{\ep\downarrow0} \frac{Ee(N_\ep\cup F_\ep)^2}
{Ee(N_\ep)^2} = 1.
\end{equation}
\citet{Ste02} argues that a useful asymptotic approach is to let
the smallest distance from the observations to the predictand
tend to 0 as $\ep\downarrow0$.
Specifically,
\citet{Ste02} proves (\ref{main}) when, essentially, for some
$x_0\in
\R^d$
not in the integer lattice, $F_\ep$ is all
points of the form $\ep(x_0+J)$ for $J$ in the integer lattice, $N_\ep$
is the restriction of $F_\ep$ to some fixed region with 0 in its interior
and $f$ is regularly varying at infinity
[\citet{BinGolTeu87}]
in every direction with a common index of variation.
The methods used in \citet{Ste02} make strong use of the gridded
nature of the observations and are not applicable here.
Furthermore, requiring
$f$ to be regularly varying at infinity with common index of variation in
all directions excludes
models for spatial-temporal phenomena that
exhibit a different degree of smoothness in space than in time.
Section \ref{sec4} provides further discussion of these issues.
Ramm (\citeyear{Ram05}), Chapter 5, takes a different approach to
studying an asymptotic
screening effect by considering a process observed with white noise
everywhere in some domain and letting the variance of the white noise
tend to 0.
In this work, we take a closer look at how the set where $Z$ is observed
affects whether an asymptotic screening effect holds.

We will take the sets $N_\ep$ and $F_\ep$ to have a particular form
that simplifies the asymptotic analysis.
Suppose $x_1,\ldots,x_n$ are distinct nonzero elements of~$\R^d$,
$y_1,\ldots,y_m$ are
distinct elements of $\R^d$ and $y_0\in\R^d$ is nonzero.
For the rest of this work, let $N_\ep= \{\ep x_1,\ldots,\ep x_n\}$
and $F_\ep=\{y_0+\ep y_1,\ldots,y_0+\ep y_m\}$.
Section \ref{sec2} explores when (\ref{main}) holds through a series of examples
leading to a~broad conjecture under
a key assumption on the spectral density $f$ of the random field:
for every $R<\infty$,
%
%
\begin{equation}\label{f-cond}
\lim_{\omega\to\infty} \sup_{|\nu|<R} \biggl| \frac{f(\omega
+\nu)}{f(\omega)}
-1\biggr| = 0.
\end{equation}
The examples will demonstrate that one generally needs a further
condition on $N_\ep$ depending on the mean square differentiability
properties of the process. For nondifferentiable processes, no further
assumptions on $N_\ep$ may be needed. Indeed, for nondifferentiable
processes on $\R$, Theorem \ref{thm1} in Section~\ref{sec3} has
(\ref{main}) as its conclusion under (\ref{f-cond}) and a mild
additional condition on $f$. For nondifferentiable processes on~$\R^2$,
if one restricts the cardinality of $N_\ep$ to~1 and of $F_\ep$ to 2
(and sets $y_2=0$), then Theorem \ref{thm2} proves (\ref{main}) under~(\ref{f-cond}) without any additional conditions on $f$.

Mat\'ern models [\citet{Ste99N1}] appear in both the examples
and the proof of Theorem \ref{thm1}.
Define $\mathcal{K}_\nu$ to be the modified Bessel function of the second
kind of order $\nu$ [\citet{autokey8}].
The Mat\'ern model on $\R^d$ has autocovariance function
$\phi(\alpha|x|)^\nu\mathcal{K}_\nu(\alpha|x|)$
for positive $\phi,\alpha$ and $\nu$.
The parameter $\nu$ controls the smoothness of the process: $Z$ has $m$
mean square derivatives in any direction if and only if $ \nu> m$.
The corresponding spectral density equals
$\phi(\alpha^2+|\omega|^2)^{-\nu-d/2}$ times a constant depending
on $\alpha,\nu$ and $d$.
All Mat\'ern models satisfy (\ref{f-cond}).

\section{Examples}\label{sec2}
This section studies a number of examples to gain some insight into
the conditions on $f$ and $N_\ep$ that are needed in order for
(\ref{main}) to hold.
The derivations of these results are elementary but not necessarily easy.
Rather than give detailed derivations of all of them, I will outline
derivations in a few of the more difficult examples in Section \ref{sec51}.

%
\begin{figure}[b]

\includegraphics{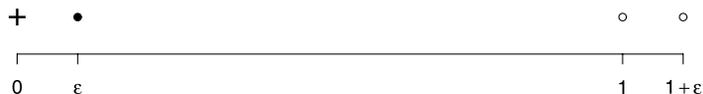}

\caption{Prediction problem for triangular autocovariance function.
Prediction site \mbox{($+$ sign)}, nearby observation (solid circle)
and distant observations (open circles).}
\label{fig1}
\end{figure}

To see why a condition like (\ref{f-cond}) is needed, let us first
consider an example on $\R$
addressed in \citet{SteHan89} and Stein (\citeyear{Ste99N1}),
pages \mbox{67--69}.
Suppose $n=1$, $x_1=1$, $m=2$, $y_0=1$, $y_1=0$ and $y_2=1$; see Figure
\ref{fig1}.
Consider $K(x)= e^{-|x|}$,
a Mat\'ern model with smoothness parameter $\frac{1}{2}$.
The corresponding process is mean square continuous but is not mean
square differentiable, and it is easy to show $Ee(N_\ep)^2\sim2\ep$
as $\ep\downarrow0$.
This process is Markov, so that $Ee(N_\ep\cup F_\ep)^2=Ee(N_\ep)^2$
for all $\ep<1$ and (\ref{main}) holds trivially.
Next consider $K(x) = (1-|x|)^+$ (where the superscipt $+$ indicates positive
part), for which $f(\omega) = \frac{1-\cos\omega}{\pi\omega^2}$,\vspace*{1pt} which
does not satisfy (\ref{f-cond}).
\citet{SteHan89}, page 180, give the BLP based on $Z(N_\ep\cup
F_\ep)$,
from which it is not difficult to show that $Ee(N_\ep)^2\sim2\ep$, just
like for $K(x)= e^{-|x|}$,
but $Ee(N_\ep\cup F_\ep)^2\sim\frac{3}{2}\ep$ as $\ep\downarrow
0$ so
that $Ee(N_\ep\cup F_\ep)^2/Ee(N_\ep)^2\to
\frac{3}{4}$ as $\ep\downarrow0$.
The choice of $y_0=1$ is critical here: for $y_0\ne1$ but positive
(keeping $x_1=1,y_1=0,y_2=1$),
$Ee(N_\ep\cup F_\ep)^2/Ee(N_\ep)^2\to1$ as $\ep\downarrow0$.
The anomaly for $y_0=1$ is related\vadjust{\goodbreak} to the lack of differentiability
of $K(x)$ at $x=1$, which is in turn related to the oscillations at
high frequencies in $f$.
See \citet{Ste05} for further discussion on the relationship of the
differentiability of $K$ away from the origin and the high-frequency behavior
of $f$.

\begin{figure}

\includegraphics{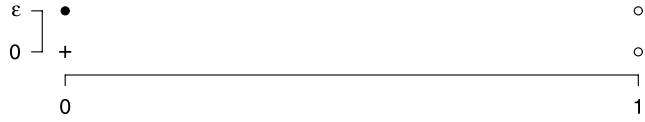}

\caption{Prediction problem for autocovariance function $K(s,t) =
e^{-|s|-|t|}$. Symbols as in Figure~\protect\ref{fig1}.}
\label{figtensor}
\end{figure}

Proposition 1 in \citet{Ste05} provides a second example showing
why a condition like (\ref{f-cond}) is needed to have a screening
effect. The following special case of this result suffices to
illustrate the point. Suppose\vspace*{1pt} $Z$ is a~stationary process
on $\R^2$ with autocovariance function $K(s,t) = e^{-|s|-|t|}$ for
$s,t\in\R$. The corresponding spectral density $f(\omega_1,\omega_2)$
is proportional to $\frac {1}{(1+\omega_1)^2 (1+\omega_2^2)}$, which
does not satisfy~(\ref{f-cond}). Consider the situation pictured in
Figure \ref{figtensor}, for which $x_1=(0,1)$, $y_0=(1,0)$, $y_1=(0,1)$
and $y_2=(0,0)$. Then using either direct calculation or Proposition 1
in \citet{Ste05}, $\lim_{\ep\downarrow0} Ee(N_\ep\cup F_\ep)^2/
Ee(N_\ep)^2 = 1-e^{-2}$.\vspace*{2pt}

%

%
\begin{figure}[b]

\includegraphics{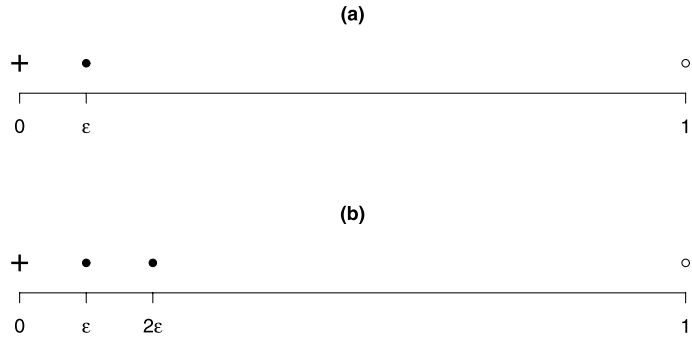}

\caption{Prediction problems for Mat\'ern model with $\nu=\frac{3}{2}$
on $\R$. Symbols as in Figure \protect\ref{fig1}.}
\label{fig2}
\end{figure}

The remaining examples all consider $f$ satisfying (\ref{f-cond}). To
see why an additional condition on $N_\ep$ is needed for (\ref{main})
to hold for differentiable processes, consider a Mat\'ern model with
smoothness parameter $\frac{3}{2}\dvtx K(x) =e^{-|x|}(1+|x|)$, for
which the corresponding process is exactly once mean square
differentiable. For $N_\ep= \{\ep\}$, $F_\ep=\{1\}$ (top
plot\vspace*{1pt} in Figure \ref{fig2}), straightforward calculations
yield $Ee(N_\ep)^2\sim \ep^2$ and $Ee(N_\ep\cup
F_\ep)^2\sim\frac{e^2-5}{e^2-4}\ep^2$ as $\ep\downarrow0$ so
$Ee(N_\ep\cup F_\ep)^2/Ee(N_\ep)^2\to \frac{e^2-5}{e^2-4}$ as
$\ep\downarrow0$. Unlike the triangular case, there is nothing special
about $y_0=1$ here and the more general result for $y_0>0$ is
$Ee(N_\ep\cup F_\ep)^2/Ee(N_\ep)^2\to1-y_0^2/(e^{2y_0}-1-2y_0-y_0^2)$.
The reason the limit is less than 1 is not because there is anything
unusual about $f$, but rather that $N_\ep$ is inadequate. Specifically,
since $Z(0) = Z(\ep) - \ep Z'(0) + o_p(\ep)$ and
$\cov\{Z(\ep),Z'(0)\}\to0$ as $\ep\downarrow0$, it is apparent that
having even a somewhat informative predictor for $Z'(0)$\vspace*{-1pt} would provide
useful information about~$Z(0)$ not contained in $Z(\ep)$. In fact, as
$\ep\downarrow0$, it is possible to show that $\hat{Z'(0)} =
\frac{e}{e^2-4}Z(\ep)-\frac{2}{e^2-4}Z(1)$\vspace*{-1pt} is an asymptotically optimal
predictor of $Z'(0)$ based on $(Z(\ep),Z(1))$ and, in turn, that
$Z(\ep)-\ep\hat{Z'(0)}$ is an asymptotically optimal predictor of
$Z(0)$ based on $(Z(\ep),Z(1))$. A screening effect does hold if $2\ep$
is added to $N_\ep$ (bottom plot of Figure \ref{fig2}). Then it is
possible to show that $Ee(N_\ep\cup F_\ep)^2\sim Ee(N_\ep)^2
\sim\frac{8}{3}\ep^3$\vspace*{2pt} as $\ep\downarrow0$, so (\ref{main}) is
true. Furthermore, as $\ep\downarrow0$, $2Z(2\ep)-Z(\ep) =
Z(\ep) - \ep[\{Z(2\ep)-Z(\ep)\}/\ep]$ is an asymptotically optimal
predictor of $Z(0)$ based on $Z(N_\ep\cup F_\ep)$ and
$\{Z(2\ep)-Z(\ep)\}/\ep$ is a~consistent predictor of $Z'(0)$. A
reasonable conjecture for a process on $\R$ with exactly~$p$ mean
square derivatives whose spectral density satisfies (\ref{f-cond}) is
that any distinct $x_1,\ldots,x_n$ with $n>p$ suffices to make
(\ref{main}) true.

It is helpful to consider this problem in the spectral domain.
We need some further notation to proceed.
For nonnegative-valued functions $a$ and~$b$ defined on a common domain
$D$, write $a(x) \ll b(x)$ if there exists finite~$C$ such that
$a(x) \le Cb(x)$ for all $x\in D$ and, for $x\in\R$,
$a(x) \ll b(x)$ as $x\downarrow x_0$ if, for some $c>0$, $a(x) \ll b(x)$
for $D=(x_0,x_0+c)$.
Write $a(x) \asymp b(x)$ if $a(x) \ll b(x)$ and $b(x) \ll a(x)$ and
define $a(x) \asymp b(x)$ as $x\downarrow0$ if $a(x) \ll b(x)$ as
$x\downarrow x_0$ and $b(x) \ll a(x)$ as $x\downarrow x_0$.
For a complex-valued function~$g$ and a nonnegative function $f$ defined
on a domain $D$ (always~$\R^d$ here),
define $\|g\|_f = \sqrt{\int_D |g(x)|^2f(x) \,dx}$.
To each random variable of the form $\sum_{j=1}^n \lambda_j Z(s_j)$
there is a corresponding
function $\sum_{j=1}^n \lambda_j e^{i\omega\cdot s_j}$,
and the mapping is an isometric isomorphism in the sense that
$E\{\sum_{j=1}^n \lambda_j Z(s_j)\}^2 = \int_{\R^d} |
{\sum_{j=1}^n \lambda_j e^{i\omega\cdot s_j}}|^2 f(\omega)
\,d\omega$.
Write $\sum_{j=1}^n \phi_{j\ep} Z(\ep x_j)$ for the BLP of $Z(0)$ based
on $Z(N_\ep)$ and
$\phi_\ep(\omega) = \sum_j \phi_{j\ep}
e^{i\ep\omega\cdot x_j}$ for the corresponding function.
If we set $\eta_\ep(\omega)=1-\phi_\ep(\omega)$, then
$Ee(N_\ep)^2= \|\eta_\ep\|^2_f$.

For any $A\subset\R^d$, call $\int_A |\eta_\ep(\omega)|^2
f(\omega) \,d\omega/\|\eta_\ep\|^2_f$ the fraction of
$Ee(N_\ep)^2$ attributable to the set of frequencies $A$.
Write $b(r)$ for the ball of radius $r$ centered at the origin.
For the scenario in Figure \ref{fig2}(a),
for any fixed $\omega_0>0$, as $\ep\downarrow0$,
%
%
\begin{equation}\label{lowfreq}
\frac{\int_{b(\omega_0)} |\eta_\ep(\omega)|^2 f(\omega)
\,d\omega}
{\|\eta_\ep\|_f^2} \sim\frac{2}{\pi}\biggl\{\tan^{-1} \omega_0
-\frac{\omega_0}
{1+\omega_0^2}\biggr\} > 0
\end{equation}
so that an asymptotically nonnegligible fraction of $Ee(N_\ep)^2$ is
attributable to a fixed range of frequencies. Similar to the definition
of $\eta_\ep$, let $\psi_\ep$ be the function corresponding to
$e(N_\ep\cup F_\ep)$, so that $Ee(N_\ep\cup F_\ep)^2 =
\|\psi_\ep\|_f^2$. Then (\ref{lowfreq}) allows $Z(1)$ to improve the
prediction nonnegligibly by making $|\psi_\ep(\omega)|^2
/|\eta_\ep(\omega)|^2$ substantially\vadjust{\goodbreak} smaller than 1 in a neighborhood
of the origin. In contrast, for the scenario in Figure \ref{fig2}(b),
$\int_{b(\omega_0)} |\eta_\ep(\omega)|^2 f(\omega) \,d\omega\ll\ep^4$
as $\ep \downarrow0$ for any fixed $\omega_0$, so that
$\int_{b(\omega_0)} |\eta_\ep(\omega)|^2 f(\omega) \,d\omega
\ll\ep\|\eta_\ep\|_f^2$ as $\ep\downarrow0$. In this\vspace*{1pt} case,
making $|\psi_\ep(\omega)|^2 /|\eta_\ep(\omega)|^2$ substantially
smaller than 1 in a~neighborhood of the origin cannot yield a~nonnegligible asymptotic impact on the mean squared prediction error.
Thus, $\int_{b(\omega_0)^c} |\psi_\ep(\omega)|^2 f(\omega)
\,d\omega/\break\int_{b(\omega_0)^c} |\eta_\ep(\omega)|^2 f(\omega)
\,d\omega$ must\vspace*{1pt} be bounded by some constant less than 1 as
$\ep\downarrow0$ for all $\omega_0$ for (\ref{main}) not to hold. The
fact that $f$ is well behaved at high frequencies [i.e., satisfies
(\ref{f-cond})] effectively precludes this possibility so that
(\ref{main}) holds. This line of reasoning forms the basis of the proof
of Theorem \ref{thm1}; see Section \ref{sec52}.

It is interesting to reconsider the two cases pictured in
Figure \ref{fig2},
for a process that is\vspace*{1pt} not quite mean square differentiable:
$K(x) = |x|\mathcal{K}_1(|x|)$, a Mat\'ern model with
smoothness parameter 1, for which $K(x) = 1 + \frac{1}{2}x^2\log(
\frac{1}{2}|x|)+\frac{1}{4}(2\gamma-1)x^2+O({x^4\log}|x|)$ as
$x\to0$
with $\gamma$ being Euler's constant.
The corresponding spectral density $f$ is proportional to
$(1+\omega^2)^{-3/2}$.
Since the process has no mean square derivatives, I conjecture that
(\ref{main}) should hold for any nonempty $N_\ep$.
For the scenario in Figure \ref{fig2}(a),
$Ee(N_\ep)^2\sim-\ep^2\log\ep$ and, for fixed
$\omega_0>0$,
\[
\int_{b(\omega_0)} |\eta_\ep(\omega)|^2 f(\omega) \,d\omega
\asymp\int_{b(\omega_0)} \frac{|1-e^{i\ep\omega}|^2 + \{K(\ep
)-1\}^2}
{(1+\omega^2)^{3/2}}\,d\omega\asymp\ep^2
\]
as $\ep\downarrow0$. Thus, the fraction of the $Ee(N_e)^2$ attributable
to $b(\omega_0)$ tends to~0 as $\ep\downarrow0$, although at only a
logarithmic rate. Not coincidentally, direct calculation shows that for
$F_\ep= \{1\}$, (\ref{main}) holds and I would expect it to hold for
more general~$F_\ep$. In fact, Theorem \ref{thm2} in Section \ref{sec3}
applies in this case and it follows that (\ref{main}) holds when
$F_\ep$ has two points (and $y_2=0$).

\begin{figure}

\includegraphics{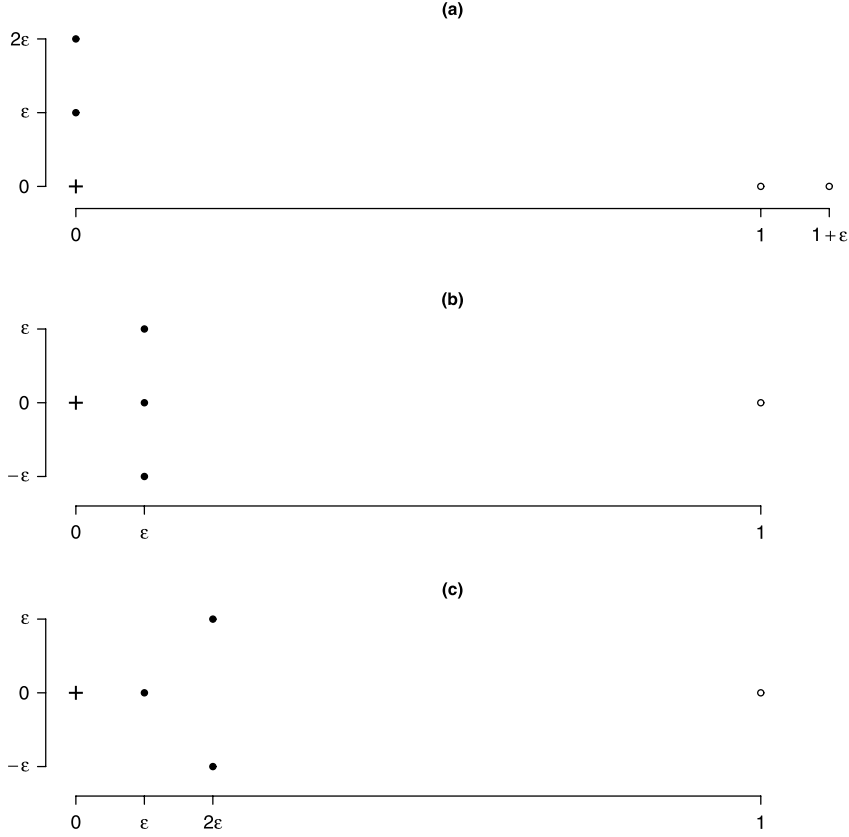}

\caption{Prediction problems for isotropic Mat\'ern model on $\R^2$ with
$\nu=\frac{3}{2}$. Symbols as in Figure~\protect\ref{fig1}.}
\label{fig3}
\end{figure}

Next consider some settings for the Mat\'ern model with
$\nu=\frac{3}{2}$ on $\R^2$. Figure \ref{fig3}(a) shows a situation in
%
%
which there are two nearby observations in the vertical direction from
the origin and two distant observations in
the horizontal direction. One might imagine that because the nearby
observations provide no information about how the process varies in the
horizontal direction, the distant observations might provide
nonneglible new information about $Z(0)$. However, Section \ref{sec51}
demonstrates that (\ref{main}) does hold in this case. The next two
examples are related to the one-dimensional examples considered in
Figure \ref{fig2} for a Mat\'ern model with $\nu=\frac{3}{2}$. Write
$Z_{i,j}$ for the $ij$th partial derivative of $Z$. In Figure
\ref{fig3}(b), $N_\ep$ has three observations, but they are collinear
along a line that does not go through the origin and it is possible to
show that the BLP of $Z_{1,0}(0,0)$ based on $Z(N_\ep)$ has
asymptotically negligible correlation with $Z_{1,0}(0,0)$ as
$\ep\downarrow0$. As a consequence, the asymptotic results are
identical to what we had in Figure \ref{fig2}(a):
$Ee(N_\ep)^2\sim\ep^2$ and $Ee(N_\ep\cup
F_\ep)^2\sim\frac{e^2-5}{e^2-4}\ep^2$ as $\ep\downarrow0$. If $N_\ep$
has three points arranged\vadjust{\goodbreak} as in Figure~\ref{fig3}(c), then
$\{Z(\ep,0)-\frac{1}{2}Z(2\ep,\ep)-\frac{1}{2}Z(2\ep,-\ep)\} /\ep$ is a
consistent predictor of $Z_{1,0}(0,0)$ and (\ref{main}) holds; see
Section \ref{sec51}.

Now consider a model satisfying (\ref{f-cond}) for which the process
is not
equally differentiable in all directions.
\citet{Ste05} gives an example of such a model.
Specifically, consider a space--time model on $\R^3\times\R$ with spectral
density $\{(1+|\omega_1|^2)^2+\omega_2^2\}^{-2}$,
$(\omega_1,\omega_2)\in\R^3\times\R$.
Writing $\erfc$ for the complementary error function,
the corresponding autocovariance function $K$ is \citet{Ste05}
%
%
\begin{eqnarray}\label{stein2005}
K(x,t)
& = & \frac{1}{16} \pi^2 e^{|x|}\erfc\biggl(|t|^{1/2}+\frac{|x|}{
2|t|^{1/2}}\biggr)\biggl(1-|x|+\frac{4t^2}{|x|}\biggr)
\nonumber\\
& &{} + \frac{1}{16} \pi^2 e^{-|x|}\erfc\biggl(|t|^{1/2}-\frac{|x|}{
2|t|^{1/2}}
\biggr)\biggl(1+|x|-\frac{4t^2}{|x|}\biggr)
\\
& &{} + \frac{1}{4}\pi^{3/2}|t|^{1/2}\exp\biggl(-|t|-\frac{|x|^2}
{4|t|}\biggr)\nonumber
\end{eqnarray}
for $x\ne0$ and $t\ne0$.
For $x=0$ or $t=0$, we can define $K$ by continuity.
For $t=0$, we get
$K(x,0)=\frac{1}{8}\pi^2e^{-|x|}(1+|x|)$, the Mat\'ern model with
$\nu=\frac{3}{2}$, so
the corresponding process is exactly once mean square differentiable in any
spatial direction.
\citet{Ste05} shows that
$K(0,t)=\frac{1}{8}\pi^2-\frac{2}{3}\pi^{3/2}|t|^{3/2}+O(t^2)$
as $t\to0$ so that $K(0,t)$ is not twice differentiable
in $t$ at $t=0$, and the corresponding process is not mean square
differentiable in time.

For (\ref{stein2005}), let us again consider the setting in Figure
\ref{fig3}(c)
with the horizontal axis corresponding to the first spatial coordinate
and the vertical axis corresponding to time.
It now turns out that the two points in $N_\ep$ off of the horizontal
axis contribute negligibly to the BLP whether or not $F_\ep$ is included.
The problem is that the lack of differentiability of $Z$ in the vertical
direction implies that the BLP of $Z_{1,0}(0,0)$ based on $Z(N_\ep)$
has asymptotic correlation 0 with $Z_{1,0}(0,0)$.
Consequently, the asymptotic results are the same as in Figure \ref
{fig2}(a) for
$K(x) = e^{-|x|}(1+|x|)$; that is,
$Ee(N_\ep\cup F_\ep)^2/Ee(N_\ep)^2\to
\frac{e^2-5}{e^2-4}$ as $\ep\downarrow0$
(Section \ref{sec51}).

\begin{figure}

\includegraphics{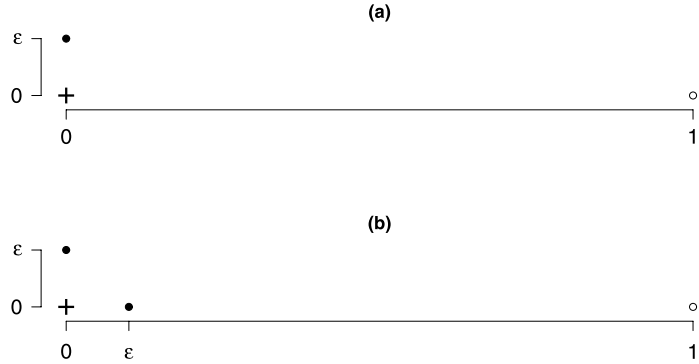}

\caption{Prediction problems for autocovariance function given in
(\protect\ref{stein2005}). Horizontal axes are first spatial coordinate
and vertical axes are time. Symbols as in Figure \protect\ref{fig1}.}
\label{fig4}
\vspace*{3pt}
\end{figure}

Figure \ref{fig4} displays two other settings we now consider for $K$
as in (\ref{stein2005}). In Figure~\ref{fig4}(a), we have $Ee(N_\ep\cup
F_\ep)^2/Ee(N_\ep)^2\to 1$ as $\ep\downarrow0$ and, in Figure~\ref{fig4}(b),
$Ee(N_\ep\cup F_\ep)^2/Ee(N_\ep)^2\to
\frac{e^2-5}{e^2-4}$ as $\ep\downarrow0$. These two cases show that it
is possible to have sets $N_\ep\subset \tilde N_\ep$ yet have that
(\ref{main}) holds for the pair of sets $(N_\ep,F_\ep)$ but not
$(\tilde N_\ep, F_\ep)$, further complicating any search for a general
result that applies to processes that are not equally smooth in all
directions.\looseness=1

%

These examples demonstrate that any general theorem that encompasses all
of them will need a condition on $N_\ep$ that depends on $f$.
The following conjecture is in accord with all of the examples
presented here:

\begin{conjecture}\label{conjec1}
Suppose $f$ satisfies (\ref{f-cond}) and the following assumption:
\renewcommand{\theassumption}{A}
\begin{assumption}\label{assumpA}
for $j=1,\ldots,n$, all
mean square derivatives of $Z$ at the origin in the direction $x_j$
can be predicted based on $Z(N_\ep)$ with mean squared error tending
to 0
as $\ep\downarrow0$.
\end{assumption}

Then for all $r>0$,
\[
\lim_{\ep\downarrow0} \frac{Ee\{N_\ep\cup b(r)^c\}^2}
{Ee\{N_\ep\}^2} = 1.
\]
\end{conjecture}

Note that here I have expanded the set of distant observations to
include all locations more than $r$ from the origin, which simplifies
the statement of the result although undoubtedly complicates its proof
(assuming it is true). It is somewhat unsatisfying to have the
condition on $N_\ep$ given in terms of properties of predictors of
derivatives of $Z$ rather than some purely geometric condition, but I
see no way to accommodate the examples treated here for $K$ as in
(\ref{stein2005}) without a condition something like Assumption
\ref{assumpA}. Verifying whether Assumption \ref{assumpA} holds in any
particular setting may require a fair amount of work, although for
$N_\ep$ of fixed and finite cardinality as we consider here, it should
generally be possible to make this determination. Note that if all mean
square derivatives of $Z$ at the origin can be consistently predicted
based on $Z(N_\ep)$ as $\ep\downarrow0$, then Assumption \ref{assumpA}
holds for any $\tilde N_\ep= \{\ep s_1,\ldots,\ep s_\ell\}$ with
$\{x_1,\ldots,x_n\}\subset\{ s_1,\ldots, s_\ell\}$.

In all of the examples for which (\ref{main}) holds,
%
%
\begin{equation}\label{low-freq}
\lim_{\ep\downarrow0}
\frac{\int_{b(\omega_0)}|\eta_\ep(\omega)|^2f(\omega) \,d\omega}
{\|\eta_\ep\|_f^2}= 0
\end{equation}
for all $\omega_0>0$, and I suspect that Assumption \ref{assumpA} is equivalent to
(\ref{low-freq}).
Examining the proof of Theorem \ref{thm1} in Section \ref{sec52}
[see (\ref{term1})], one sees that (\ref{low-freq})
is essential to making the proof work.

\section{Theorems}\label{sec3}

I do not know how to prove Conjecture \ref{conjec1} in anything like its full generality.
Assuming $Z$ is not differentiable in any direction simplifies matters
considerably, because Assumption \ref{assumpA} then holds for any~nonemp\-ty~$N_\ep$.
Theorem~\ref{thm1} considers nondifferentiable processes on $\R$
and Theorem~\ref{thm2} nondifferentiable processes on $\R^2$.

\begin{theorem}\label{thm1}
Suppose, for $d=1$ and some $\alpha\in(0,2)$,
%
%
\begin{equation}\label{fasym}
f(\omega) \asymp(1+|\omega|)^{-\alpha-1},
\end{equation}
and $f$ satisfies (\ref{f-cond}).
Then (\ref{main}) holds.
\end{theorem}

Condition (\ref{fasym}) is stronger than necessary to guarantee $Z$
is not differentiable.
Because part of the proof is to show that the low frequencies do not
matter in the limit, (\ref{fasym}) can likely be weakened
to hold only for all $\omega$ sufficiently large.
Removing~(\ref{fasym}) entirely would be more difficult.\vadjust{\goodbreak}

The next theorem applies to nondifferentiable processes
in $\R^2$ and does not require any conditions on $f$ beyond (\ref{f-cond}).
However, it does restrict $N_\ep$ to have only one point and $F_\ep$
to have
two.
The theorem also assumes $y_2=0$, but this restriction does not
meaningfully detract from the content of the result and, in any case,
could be removed at the cost of a somewhat messier proof.
Extending the result to $\R^d$ is straightforward,
but taking $d>3$ is pointless
in this setting because any 4 points in $\R^d$ fall on a three-dimensional
hyperplane, and even taking $d=3$ provides no new insight beyond what is
learned from the two-dimensional setting.
\begin{theorem}\label{thm2}
Suppose $Z$ has spectral density $f$ satisfying (\ref{f-cond}) and
that~$Z$ is not mean square differentiable in any direction.
In addition, suppose $N_\ep= \{\ep x_1\}$ and $F_\ep=\{y_0,y_0+\ep
y_1\}$,
where $x_1,y_0$ and $y_1$ are all nonzero.
Then (\ref{main}) holds.
\end{theorem}

Note that the example
referred to in Figure \ref{figtensor} satisfies the conditions
on~$N_\ep$ and $F_\ep$ in Theorem \ref{thm2}, and the process is not mean square
differentiable in any direction, but $f$ does not satisfy (\ref{f-cond}).
As we have seen,
(\ref{main}) does not hold in this setting, so that Theorem \ref{thm2} would be false
if we removed (\ref{f-cond}).

Throughout this work we assume that $Z$ has a known mean 0.
It is common in practice to assume that $Z$ has an unknown constant
mean $\mu$ and then predict $Z(0)$ by what is called the ordinary
kriging predictor, which is just an example of
the best linear unbiased predictor [Stein (\citeyear{Ste99N1})].
In all of the examples considered in Section \ref{sec2}, for which (\ref{main}) holds
for simple kriging, it still holds for ordinary kriging.
Furthermore, Theorems \ref{thm1} and \ref{thm2} can be easily shown to hold for ordinary kriging
by proving that, under the conditions of the theorems, the ordinary
kriging predictor based on $N_\ep$ is asymptotically optimal relative
to the simple kriging predictor (see the ends of each proof in Section \ref{sec5}).
Thus, if Conjecture \ref{conjec1} holds for simple kriging, then I would expect it
also holds for ordinary kriging.

\section{Discussion}\label{sec4}

The space--time process on $\R^3\times\R$ considered in Section
\ref{sec3}
with spectral density $\{(1+|\omega_1|^2)^2+\omega_2^2\}^{-2}$,
$(\omega_1,\omega_2)\in\R^3\times\R$, is an example of a process
with a different degree of differentiability in time than in space.
It is a special case of the stochastic fractional heat equations
studied by \citet{KelLeoRui05}, which are in
turn a special case of a class of space--time processes suggested in
\citet{Ste05} whose spectral densities are of the form
%
%
\begin{equation}\label{doubly-matern}
f(\omega_1,\omega_2) = \{c_1(a_1^2+|\omega_1|^2)^{\alpha_1}+
c_2(a_2^2+|\omega_2|^2)^{\alpha_2}\}^{-\nu}
\end{equation}
for\vspace*{1pt} $\omega_1\in\R^{d_1}$, $\omega_2\in\R^{d_2}$, $\nu>
\frac {d_1}{2\alpha_1} +\frac{d_2}{2\alpha_2}$ and
$c_1,c_2,\alpha_1,\alpha_2$ and $a_1^2+a_2^2$ positive to ensure $f$ is
integrable. Because of the superficial similarity of this model to the
Mat\'ern model, we might call it doubly Mat\'ern. All spectral
densities of the form (\ref{doubly-matern}) satisfy (\ref{f-cond}) and
thus, I conjecture, satisfy an asymptotic screening effect whenever
Assumption \ref{assumpA} applies to $N_\ep$. At the same time, by
adjusting the parameters $\alpha_1,\alpha_2$ and $\nu$, we can obtain
processes with any desired degree of differentiability in time and any
separate degree of differentiability in space [\citet{Ste05}].
Note that $f$ of the form (\ref{doubly-matern}) satisfies the
conditions of Theorem \ref{thm2} when $d_1=d_2=1$,
$2\nu\le\frac{3}{\alpha_1}+\frac {1}{\alpha_2}$ and
$2\nu\le\frac{1}{\alpha_1}+\frac{3}{\alpha_2}$, the last two conditions
being necessary and sufficient to make $Z$ not mean square
differentiable in any direction. \citet{autokey15} derives some
results for the covariance structure when $a_1=a_2=0$ and $\alpha_2=1$.

Despite its flexibility, model (\ref{doubly-matern}) is still
restrictive in some ways, in
particular in exhibiting what \citet{Gne02} calls full symmetry,
due to the fact that $f(\omega_1,\omega_2)=f(\omega_1,-\omega_2)$, and
hence the corresponding process has the same covariance structure
with time running backwards as it does with time running forward.
Thus, for example, this model is unsuitable for processes with a
dominant direction of advection.
\citet{Ste05} discusses possible approaches to extending this
model to
allow for asymmetries.

As noted in Section \ref{sec3}, (\ref{low-freq}), which says that only an
asymptotically negligible fraction of $Ee(N_\ep)^2$
can be attributed to some fixed frequency range, is crucial to obtaining
a screening effect.
This same property was also the key idea in \citet{Ste99N2} to obtaining
explicit results on the asymptotic efficiency of predictors based on
an incorrect spectral density having similar behavior to the correct spectral
density at high frequencies.
The high-frequency behavior of a Gaussian process is also crucial
to estimation of the covariance structure [\citet{Ste99N1}], and
misspecification
of this high-frequency behavior can lead to poor behavior of estimates,
particularly if likelihood-based methods are used [\citet{Ste99N1},
Chapter~6,
and \citet{Ste08}].
As statisticians strive to advance the statistical analysis of
spatial-temporal processes, they should pay close attention to the
spectral behavior of the models they use.
In particular, models that do not satisfy (\ref{f-cond}) should be
used with caution.

\section{Proofs}\label{sec5}

\subsection{Examples}\label{sec51}

For a random vector $Y$, write $\cov(Y)$ for the covariance matrix of $Y$,
write 0 for a column vector of zeroes whose length is apparent from
context and denote transposes by primes.
The following result simplifies the calculations for several of the examples.
\begin{lemma}\label{lem1}
$\!\!\!$If there exists $a(\ep)\,{>}\,0$, $\delta_\ep\,{\in}\,\R^{n+m}$ and $\Delta
_\ep$ an
\mbox{$(n\,{+}\,m)\,{\times}\,(n\,{+}\,m)$} matrix such that
%
%
\begin{equation}\label{lemma-limit}
\lim_{\ep\downarrow0}
\cov\pmatrix{
a(\ep)\{Z(0) - \delta_\ep\cdot Z(N_\ep\cup F_\ep)\} \cr
\Delta_\ep Z(N_\ep\cup F_\ep)}
= \pmatrix{
k & 0' \cr
0 & K
}
\end{equation}
for some $k>0$ and $K$ positive definite, then
\[
\lim_{\ep\downarrow0}
\frac{E\{Z(0) - \delta_\ep\cdot Z(\ep)\}^2}{Ee(N_\ep\cup F_\ep
)^2} = 1.
\]
\end{lemma}

For (\ref{lemma-limit}) to hold, $\tilde e_\ep= Z(0) - \delta_\ep\cdot
Z(N_\ep\cup F_\ep)$ must satisfy $E e(N_\ep\cup F_\ep)^2/\break E\tilde
e_\ep^2 \to1$ as $\ep\downarrow0$. To prove the lemma, note that
(\ref{lemma-limit}) and $K$ positive definite imply $\cov\{ \Delta_\ep
Z(\ep)\}$ is positive definite for all $\ep$ sufficiently small. Thus,
for all~$\ep$ sufficiently small, the BLP of $Z(0)$ based on
$Z(N_\ep\cup F_\ep)$ is the same as the BLP of $Z(0)$ based on
$\Delta_\ep Z(N_\ep\cup F_\ep)$. Since matrix inverse is a continuous
function in some neighborhood of $K$, using basic results on BLPs
[e.g., Stein (\citeyear{Ste99N1}), Section~1.2],
\begin{eqnarray*}
& & a(\ep)^2 Ee(N_\ep\cup F_\ep)^2 \\
& &\qquad = a(\ep)^2\bigl(\var\tilde e_\ep
-\cov\{\tilde e_\ep,Z(N_\ep\cup F_\ep)'\Delta'_\ep\} \\
& &\qquad\quad\hspace*{37pt}\hspace*{25.6pt}{} \times
[\cov\{\Delta_\ep Z(N_\ep\cup F_\ep)\}]^{-1}
\cov\{\Delta_\ep Z(N_\ep\cup F_\ep),\tilde e_\ep
\}\bigr) \\
& &\qquad = \var\{a(\ep)\tilde e_\ep\}
-\cov\{a(\ep)\tilde e_\ep,Z(N_\ep\cup F_\ep)'\Delta'_\ep\} \\
& &\qquad\quad{} \times
[\cov\{\Delta_\ep Z(N_\ep\cup F_\ep)\}]^{-1}
\cov\{\Delta_\ep Z(N_\ep\cup F_\ep),a(\ep)\tilde e_\ep
\} \\
& &\qquad \to k - 0' K^{-1} 0
\end{eqnarray*}
as $\ep\downarrow0$, and the lemma follows.

To apply Lemma \ref{lem1} to the setting in Figure \ref{fig3}(a) with
$K(x)=e^{-|x|}(1+|x|)$, it suffices to show
\begin{eqnarray*}
& & \lim_{\ep\downarrow0}
\cov\pmatrix{
\ep^{-3/2}\{Z(0,0) - 2 Z(0,\ep)+Z(0,2\ep)\} \cr
Z(0,\ep) \cr
\ep^{-1}\{Z(0,2\ep)-Z(0,\ep)\} \cr
Z(1,0) \cr
\ep^{-1}\{Z(1+\ep,0)-Z(1,0)\}}
\\[3pt]
& &\qquad = \pmatrix{
\frac{8}{3} & 0 & 0 & 0 & 0 \cr
0 & 1 & 0 & 2e^{-1} & -e^{-1} \cr
0 & 0 & 1 & 0 & 0 \cr
0 & 2e^{-1} & 0 & 1 & 0 \cr
0 & -e^{-1} & 0 & 0 & 1}.
\end{eqnarray*}
To show, for example, that $\cov[ \ep^{-1}\{Z(0,2\ep)-Z(0,\ep)\},
\ep^{-1}\{Z(1+\ep,0)-Z(1,0)\}]\to0$, define the function $\tilde K$
on $[0,\infty)$ by $\tilde K(r) = e^{-r}(1+r)$, which has bounded derivatives
of all orders on $[0,\infty)$.
Then using a Taylor series,
\begin{eqnarray*}
\hspace*{-2pt}& & \cov\{Z(0,2\ep)-Z(0,\ep),Z(1+\ep,0)-Z(1,0)\} \\
\hspace*{-2pt}& &\quad = K\bigl(\sqrt{(1+\ep)^2 +4\ep^2}\bigr) - K\bigl(\sqrt
{(1+\ep)^2
+\ep^2}\bigr) \\
\hspace*{-2pt}& &\qquad{} - K\bigl(\sqrt{1+4\ep^2}\bigr) + K\bigl(\sqrt{1+\ep
^2}\bigr) \\
\hspace*{-2pt}& &\quad = K'(1+\ep)\bigl\{\sqrt{(1+\ep)^2 + 4\ep^2}-
\sqrt{(1+\ep)^2 + \ep^2}\bigr\} \\
\hspace*{-2pt}& &\qquad{} -K'(1)\bigl\{\sqrt{1+4\ep^2}-\sqrt{1+\ep^2}\bigr\} +
O(\ep^4) \\
\hspace*{-2pt}& &\quad = K'(1+\ep)(1+\ep)\biggl\{\frac{2\ep^2}{(1+\ep)^2}-\frac{\ep^2}
{2(1+\ep)^2}\biggr\} - K'(1)\biggl(2\ep^2-\frac{1}{2}\ep^2\biggr) +
O(\ep^4) \\
\hspace*{-2pt}& &\quad \ll\ep^3,
\end{eqnarray*}
and $\cov[ \ep^{-1}\{Z(0,2\ep)-Z(0,\ep)\},
\ep^{-1}\{Z(1+\ep,0)-Z(1,0)\}]\to0$ follows.

Lemma \ref{lem1} can be applied to the setting in Figure \ref{fig3}(c)
with $K(x)=e^{-|x|}(1+|x|)$ by showing
\begin{eqnarray*}
& & \lim_{\ep\downarrow0}
\cov\pmatrix{
\ep^{-3/2}\bigl\{Z(0,0) - 2
Z(\ep,0)+\frac{1}{2}Z(2\ep,\ep)+\frac{1}{2}Z(2\ep,-\ep)\bigr\} \vspace*{2pt}\cr
Z(\ep,0) \cr
\ep^{-1}\{2Z(\ep,0)-Z(2\ep,\ep)-Z(2\ep,-\ep)\} \cr
Z(1,0)}
\\
& &\qquad = \pmatrix{
\frac{1}{3}\bigl(10\sqrt{5}-8\sqrt{2}\bigr) & 0 & 0 & 0 \cr
0 & 1 & 0 & 2e^{-1} \cr
0 & 0 & 4 & -2e^{-1} \cr
0 & 2e^{-1} & -2e^{-1} & 1}.
\end{eqnarray*}
Specifically, it is not necessary to consider $Z(2\ep,\ep)$
and $Z(2\ep,-\ep)$ separately: by symmetry,
the BLP of $Z(0,0)$ based on $Z(N_\ep\cup F_\ep)$
depends on $Z(2\ep,\ep)$ and $Z(2\ep,-\ep)$ only
through $Z(2\ep,\ep)+Z(2\ep,-\ep)$.

As a final example, let us apply Lemma \ref{lem1} to the setting in Figure
\ref{fig3}(c) with $K$ given by (\ref{stein2005}). Again by symmetry, we
can restrict to predictors that depend on $Z(2\ep,\ep)$ and
$Z(2\ep,-\ep)$ only through $Z(2\ep,\ep)+Z(2\ep,-\ep)$. For $a$ and $b$
fixed and positive, using a Taylor series and
\[
\erfc(x) = 1 -\frac{2}{\sqrt{\pi}}\biggl(x- \frac{1}{3}x^3\biggr)
+ O(|x|^5)
\]
as $x\to0$ [\citet{autokey8}, page 162], it is possible to show
\[
K(a\ep,b\ep) = \tfrac{1}{8}\pi^2 - \tfrac{2}{3}(\pi b\ep)^{3/2} +
O(\ep^2)
\]
as $\ep\downarrow0$.
This result also holds when $a$ or $b$ equals 0.
It follows that
\begin{eqnarray*}
& & \lim_{\ep\downarrow0}
\cov\pmatrix{
\ep^{-1}\{Z(0,0) - Z(\ep,0)\} \cr
Z(\ep,0) \cr
\ep^{-3/4}\{2Z(\ep,0)-Z(2\ep,\ep)-Z(2\ep,-\ep)\}}
\\[3pt]
& &\qquad = \pmatrix{
\frac{1}{8}\pi^2 & 0 & 0 \vspace*{2pt}\cr
0 & \frac{1}{8}\pi^2 & 0 \vspace*{2pt}\cr
0 & 0 & \frac{8}{3}\bigl(2-\sqrt{2}\bigr)\pi^{3/2}}
\end{eqnarray*}
so that $Z(\ep,0)$ is an asymptotically optimal predictor of $Z(0,0)$
based on~$N_\ep$.
Furthermore, for $c_1=2/(e^2-4)$ and $c_2=-e/(e^2-4)$,
\begin{eqnarray*}
& & \lim_{\ep\downarrow0}
\cov\pmatrix{
\ep^{-1}\{Z(0,0) - (1+c_1\ep)Z(\ep,0)-c_2\ep Z(1,0)\} \cr
Z(\ep,0) \cr
Z(1,0) \cr
\ep^{-3/4}\{2Z(\ep,0)-Z(2\ep,\ep)-Z(2\ep,-\ep)\}}
\\[3pt]
& &\qquad = \pmatrix{
\dfrac{1}{8}\pi^2\dfrac{e^2-5}{e^2-4} & 0 & 0 & 0 \vspace*{2pt}\cr
0 & \dfrac{1}{8}\pi^2 & \dfrac{\pi^2}{4e} & 0 \vspace*{2pt}\cr
0 & \dfrac{\pi^2}{4e} & \dfrac{1}{8}\pi^2 & 0 \vspace*{2pt}\cr
0 & 0 & 0 & \dfrac{8}{3}\bigl(2-\sqrt{2}\bigr)\pi^{3/2}},
\end{eqnarray*}
and the conditions of Lemma \ref{lem1} are satisfied.

\subsection{\texorpdfstring{Proof of Theorem \protect\ref{thm1}}{Proof of Theorem 1}}\label{sec52}

Theorem 3.1 in \citet{Xia08} implies
%
%
\begin{equation}\label{rate}
\|\eta_\ep\|_f^2=Ee(N_\ep)^2 \asymp\ep^\alpha
\end{equation}
as $\ep\downarrow0$.
Let us use (\ref{rate}) to show that ${\sum_{j=1}^n} |\phi_{j\ep}|$
is bounded in
$\ep$ as $\ep\downarrow0$.
If we define $M_\ep= \max(1,\sum_{j=1}^n |\phi_{j\ep}|)$
and $x_0=0$, we can
write $\eta_\ep(\omega)$ in the form $M_\ep\sum_{j=0}^n \mu_{j\ep
} e^{i\ep\omega
x_j}$ for appropriate $\mu_{j\ep}$'s, where, by construction,
$|\mu_{j\ep}|\le1$ for all $j$ and $\ep$.
Thus, if we can show $M_\ep$ bounded, then
${\sum_{j=1}^n }|\phi_{j\ep}|$ is also bounded.
By (\ref{fasym}), there exists $0<C_1\le C_2<\infty$ such that
%
%
\begin{equation}\label{f-bounds}
\frac{C_1}{(1+|\omega|)^{\alpha+1}}\le f(\omega) \le\frac
{C_2}{(1+|\omega|)^{\alpha+1}}
\end{equation}
for all $\omega$.
Thus, making the change of variables $\nu= \ep\omega$ in the second
step,
%
%
\begin{eqnarray}\label{bound1}
Ee(N_\ep)^2 & \ge&
C_1M_\ep^2 \int_{-\infty}^{\infty}\Biggl| \sum_{j=0}^n \mu_{j\ep}
e^{i\ep\omega x_j}\Biggr|^2 (1+|\omega|)^{-\alpha-1} \,d\omega
\nonumber\\
& = & C_1M_\ep^2 \ep^\alpha\int_{-\infty}^{\infty}\Biggl| \sum_{j=0}^n
\mu_{j\ep} e^{i\nu x_j}\Biggr|^2 (\ep+|\nu|)^{-\alpha-1} \,d\nu
\\
& \ge& C_1M_\ep^2 \biggl(\frac{1}{2}\ep\biggr)^\alpha\int_1^\infty
\Biggl| \sum_{j=0}^n \mu_{j\ep} e^{i\nu x_j}\Biggr|^2 \nu
^{-\alpha-1} \,d\nu\nonumber
\end{eqnarray}
for all $\ep<1$.
Suppose $M_\ep$ is unbounded.
Then there exists a sequence $\{\ep(k)\}$ tending to 0
such that $M_{\ep(k)}\to\infty$.
Because the $\mu_{j\ep}$'s are bounded,
there exists $ (\mu_0,\ldots,\mu_n)\in\R^{n+1}$ and
a subsequence of $\{\ep(k)\}$, call it $\{\ep(k_\ell)\}$, along
which\vadjust{\goodbreak}
$(\mu_{0\ep(k_\ell)},\ldots,\mu_{n\ep(k_\ell)}) \to(\mu
_0,\ldots,\mu_n)$
as $\ell\to\infty$.
Since $\alpha>0$, by dominated convergence, it follows that
\[
\int_1^\infty
\Biggl| \sum_{j=0}^n \mu_{j\ep(k_\ell)} e^{i\nu x_j}\Biggr|^2 \nu
^{-\alpha-1}
\,d\nu\to\int_1^\infty
\Biggl| \sum_{j=0}^n \mu_{j} e^{i\nu x_j}\Biggr|^2 \nu^{-\alpha-1}
\,d\nu> 0
\]
as $\ell\to\infty$, which, together with (\ref{bound1}),
contradicts (\ref{rate}), so $M_\ep$ and
${\sum_{j=1}^n }|\phi_{j\ep}|$ must be bounded as $\ep\downarrow0$.

Now consider the behavior of $\eta_\ep$ at low frequencies.
Define $p_{\ep} = 1/\sum_j \phi_{j\ep}$ and
$\tilde\eta_\ep(\omega) = 1 - p_{\ep}\sum_{j=1}^n
\phi_{j\ep}e^{i\ep\omega x_j}$.
By (\ref{fasym}) and (\ref{rate}), $\int_0^1 |\eta_\ep(\omega)|^2
\,d\omega\ll\ep^\alpha$ and, writing $\mathrm{Re}$ for real part,
$|\eta_\ep(\omega)|^2 \ge\{\operatorname{Re} \eta_\ep(\omega)\}^2 =
(1-p_\ep)^2 +O(\ep^2)$ uniformly for $\omega\in[0,1]$.
It follows that
%
%
\begin{equation}\label{Phirate}
(p_{\ep} - 1)^2 \ll\ep^\alpha
\end{equation}
as $\ep\downarrow0$.
Using $|e^{ix}-1| \le|x|$ for all $x\in\R$,
for $\beta\in[0,1]$ and $\alpha\in(0,2)$,
%
%
\begin{equation}\label{phibound}
\int_{b(\ep^{-\beta})} | \tilde\eta_\ep(\omega)|^2 f(\omega)
\,d\omega\ll\ep^2\int_0^{\ep^{-\beta}} \frac{\omega
^2}{1+\omega^{\alpha+1}}
\,d\omega\ll\ep^{2-\beta(2-\alpha)}
\end{equation}
as $\ep\downarrow0$.
Because $\sum_{j=1}^n \phi_{j\ep} Z(\ep x_j)$ is the BLP of $Z(0)$,
$\|\eta_\ep\|^2_f\le\|\tilde\eta_\ep\|^2_f$, so that
%
%
\begin{eqnarray}\label{lowfreqbd}
\int_{b(\ep^{-\beta})} | \eta_\ep(\omega)|^2 f(\omega)
\,d\omega& \le& \int_{b(\ep^{-\beta})} | \tilde\eta_\ep(\omega)|^2
f(\omega) \,d\omega\nonumber\\[-8pt]\\[-8pt]
& &{} + \int_{b(\ep^{-\beta})^c}\{
| \tilde\eta_\ep(\omega)|^2 - | \eta_\ep(\omega)|^2\}
f(\omega) \,d\omega.\nonumber
\end{eqnarray}
Straightforward algebra shows
%
%
\begin{eqnarray}\label{eta-bound}
& & | \tilde\eta_\ep(\omega)|^2 - | \eta_\ep(\omega)|^2
\nonumber\\
& &\qquad = (p_\ep^2-1)|\phi_\ep(\omega)|^2-2(p_\ep-1)\operatorname{Re}
\phi_\ep(\omega) \\
& &\qquad = 2(p_\ep-1)^2|\phi_\ep(\omega)|^2+2(p_\ep-1)[
|\phi_\ep(\omega)|^2-\operatorname{Re} \phi_\ep(\omega)].\nonumber
\end{eqnarray}
The boundedness of the $\phi_{j\ep}$'s in $\ep$ implies
$| |\phi_\ep(\omega)|^2-p_\ep^{-2}|^2 \ll\min(1,\ep
^2\omega^2)$
and $|{\operatorname{Re} \phi_\ep}(\omega) - p_\ep^{-1}|\ll
\min(1,\ep^2\omega^2)$, and it follows that
\[
\bigl| |\phi_\ep(\omega)|^2-\operatorname{Re} \phi_\ep(\omega)\bigr|
\ll|p_\ep-1| + \min(1,\ep^2\omega^2)
\]
as $\ep\downarrow0$,
which, together with (\ref{Phirate}) and (\ref{eta-bound}), yields
\[
| \tilde\eta_\ep(\omega)|^2 - | \eta_\ep(\omega)|^2 \ll\ep
^\alpha+
\ep^{\alpha/2}\min(1,\ep^2\omega^2)
\]
as $\ep\downarrow0$. Thus,
%
%
\begin{eqnarray}\label{eta-tilde}
& & \int_{b(\ep^{-\beta})^c} \{| \tilde\eta_\ep(\omega)|^2 -
| \eta_\ep(\omega)|^2 \}f(\omega) \,d\omega\nonumber\\
& &\qquad \ll
\int_{\ep^{-\beta}}^{\ep^{-1}}
\frac{\ep^\alpha+\ep^{2+\alpha/2}\omega^2}{\omega^{\alpha
+1}}\,d\omega+
\int_{\ep^{-1}}^\infty\frac{\ep^{\alpha/2}}{\omega^{\alpha+1}}
\,d\omega
\\
& &\qquad \ll\ep^{3\alpha/2} + \ep^{\alpha(\beta+1)}\nonumber
\end{eqnarray}
as $\ep\downarrow0$.
Combining this bound with (\ref{phibound}) and
(\ref{lowfreqbd}) implies that for all $\beta\in[0,1]$,
%
%
\begin{equation}\label{lowfreqbd2}
\int_{b(\ep^{-\beta})} | \eta_\ep(\omega)|^2 f(\omega)
\,d\omega\ll\ep^{2-\beta(2-\alpha)} + \ep^{3\alpha/2} +
\ep^{\alpha(\beta+1)}
\end{equation}
as $\ep\downarrow0$.
Note that the bound in (\ref{lowfreqbd2}) is
$o(\ep^\alpha)$ as $\ep\downarrow0$
for all $\alpha\in(0,2)$ and $\beta\in(0,1)$.

Let $\Lambda_\ep=
(\la_{1\ep},\ldots,\la_{m\ep})$, and
assume $\Lambda_\ep\ne0$ hereafter, as the case $\Lambda_\ep=0$ is
trivial to
handle.
We next show
the correlation of $e(N_\ep)$ and $\Lambda_\ep\cdot Z(F_\ep)$ is
asymptotically negligible.
Defining $\la_\ep(\omega) = \sum_{j=1}^m \la_{j\ep} e^{i\ep
\omega y_j}$,
%
%
\begin{eqnarray}\label{corr1}
& &
\corr\{ e(N_\ep), \Lambda_\ep\cdot Z(F_\ep)
\} \nonumber\\
& &\qquad = \frac{\int_{b(\ep^{-\beta})} \eta_\ep(\omega) e^{-i\omega y_0}
\overline{\la_{\ep} (\omega)}f(\omega) \,d\omega}{\|\eta_\ep\|_f
\|\la_\ep\|_f}\nonumber\\[-8pt]\\[-8pt]
&&\qquad\quad{}+ \frac{\int_{b(\ep^{-\beta})^c} \eta_\ep(\omega) e^{-i\omega y_0}
\overline{\la_{\ep} (\omega)}f(\omega) \,d\omega}{\|\eta_\ep\|_f
\|\la_\ep\|_f} \nonumber\\
& &\qquad \stackrel{\Delta}{=} I_1+I_2.\nonumber
\end{eqnarray}
Using the Cauchy--Schwarz inequality and (\ref{lowfreqbd2}), for
$\beta\in(0,1)$,
%
%
\begin{equation}\label{term1}
I_1 \le
\frac{\sqrt{\int_{b(\ep^{-\beta})}|\eta_\ep(\omega)|^2
f(\omega) \,d\omega}}{\|\eta_\ep\|_f}\to0
\end{equation}
as $\ep\downarrow0$, uniformly in $\Lambda_\ep$.
Next,
define $R_k = 2\pi k/y_0$ and $k_\ep= \lfloor y_0\ep^{-\beta}/(2\pi
)\rfloor$.
Then $R_{k_\ep}\le\ep^{-\beta}$ and
%
%
\begin{eqnarray}\label{corr2}
& & \biggl|\int_{b(\ep^{-\beta})^c} \eta_\ep(\omega) e^{-i\omega y_0}
\overline{\la_{\ep} (\omega)}f(\omega) \,d\omega\biggr| \nonumber
\\
& &\qquad \le2\sum_{k=k_\ep}^\infty\biggl|\int_{R_k}^{R_{k+1}}
\eta_\ep(\omega) e^{-i\omega y_0}
\overline{\la_{\ep} (\omega)}f(\omega) \,d\omega\biggr| \nonumber
\\
& &\qquad \le2\sum_{k=k_\ep}^\infty f(R_k)\biggl|\int_{R_k}^{R_{k+1}}
\eta_\ep(\omega) e^{-i\omega y_0}
\overline{\la_{\ep} (\omega)} \,d\omega\biggr|
\\
& &\qquad\quad{} + 2\sum_{k=k_\ep}^\infty\int_{R_k}^{R_{k+1}}
|\eta_\ep(\omega)
\la_{\ep} (\omega)|
|f(\omega)-f(R_k)| \,d\omega
\nonumber\\
& &\qquad \stackrel{\Delta}{=} I_3+I_4.\nonumber
\end{eqnarray}
For $\omega\in(R_k,R_{k+1}]$, by (\ref{f-cond}),
there exist constants $c_k\to0$ as $k\to\infty$ such that
%
%
\begin{equation}\label{fbound}
|f(\omega)-f(R_k)| \le
c_k \min\{f(R_k),f(\omega)\},
\end{equation}
so
%
%
\begin{eqnarray}\label{corr3}
I_4 & \le& 2\sum_{k=k_\ep}^\infty c_k \int_{R_k}^{R_{k+1}} |\eta
_\ep(\omega)
\la_{\ep} (\omega)| f(\omega) \,d\omega\nonumber\\
& \le& 2 \sup_{k\ge k_\ep} c_k \int_{R_{k_\ep}}^\infty|\eta_\ep
(\omega)
\la_{\ep} (\omega)| f(\omega) \,d\omega\\
& \le& \sup_{k\ge k_\ep} c_k \|\eta_\ep\|_f
\|\la_\ep\|_f,\nonumber
\end{eqnarray}
the last step by the Cauchy--Schwarz inequality.
Now consider $I_3$ in (\ref{corr2}).
Defining $\theta_\ep(\omega) = \eta_\ep(\omega) \overline{\la
_\ep(\omega)}$,
we can write $\theta_\ep$ in the form
$\sum_{j=1}^{m(n+1)}\theta_{j\ep}e^{i\ep\omega z_j}$.
For $M_\ep= \max(1,{\sum_{j=1}^n} |\phi_{j\ep}|)$,
let $M = \limsup_{\ep\downarrow0} M_\ep$, which we showed is finite.
Then, setting $L_\ep={\sum_{j=1}^m }| \la_{j\ep}|$,
it is easy to show that ${\sum_{j=1}^{m(n+1)}}|\theta_{j\ep}| \le
(2M+1)L_\ep$ for all $\ep$ sufficiently small.
Integrating by parts,
\begin{eqnarray*}
\int_{R_k}^{R_{k+1}}\theta_\ep(\omega)e^{-i\omega y_0} \,d\omega
&=& \frac{e^{-iR_k y_0}}{iy_0}\{\theta_\ep(R_k)-\theta_\ep
(R_{k+1})\}\\
&&{}+\frac{1}{iy_0}\int_{R_k}^{R_{k+1}} e^{-i\omega y_0} \theta_\ep
'(\omega)
\,d\omega.
\end{eqnarray*}
Defining $\check{z} = {\max_j} |z_j|$, we have
$|\theta_\ep(R_k)-\theta_\ep(R_{k+1})| \le{\sum_j }|\theta_{j\ep}|
|1-e^{i2\pi\ep z_j/y_0}|\le2\pi(2M+1)\check{z}L_\ep\ep/y_0$
and $|\theta_\ep'(\omega)|\le
(2M+1)\check{z}L_\ep\ep$ for all $\ep$ sufficiently small, so that
%
%
\begin{equation}\label{intbound}
\biggl|\int_{R_k}^{R_{k+1}}\theta_\ep(\omega)e^{-i\omega y_0}
\,d\omega\biggr|
\le\frac{4\pi}{y_0^2}(2M+1)\check{z}L_\ep\ep
\end{equation}
for all $\ep$ sufficiently small.
Setting $\beta=\frac{1}{2}$,
inequalities (\ref{f-bounds}) and (\ref{intbound}) imply
%
%
\begin{eqnarray}\label{bound2}
I_3 & \le& 2C_2(2M+1)\check{z}L_\ep\sum_{k=k_\ep}^\infty\frac
{\ep}{k^{\alpha+1}}
\nonumber\\[-8pt]\\[-8pt]
& \le& 2\alpha^{-1}C_2(2M+1)\check{z}L_\ep\biggl(\frac
{8}{y_0}\biggr)^\alpha
\ep^{\alpha/2+1}
\nonumber
\end{eqnarray}
for all $\ep$ sufficiently small.
Similarly to (\ref{rate}), it is possible to show
$L_\ep\ep^{\alpha/2}\ll\|\la_\ep\|_f$ as $\ep\downarrow0$, so that
by (\ref{rate}) and (\ref{bound2}),
%
%
\begin{equation}\label{corr4}
\frac{I_3}
{\|\eta_\ep\|_f\|\la_\ep\|_f} \ll\ep^{1-\alpha/2}
\end{equation}
as $\ep\downarrow0$
uniformly in $\Lambda_\ep$.
Since $\alpha<2$, this bound tends to 0 uniformly in~$\Lambda_\ep$.
Applying (\ref{corr3}) and (\ref{corr4}) to (\ref{corr2})
yields $I_2$ [defined in (\ref{corr1})] tending to 0
as $\ep\downarrow0$ uniformly in~$\Lambda_\ep$, which
together with (\ref{corr1}) and (\ref{term1}), implies
%
%
\begin{equation}\label{corr5}
\lim_{\ep\downarrow0} \sup_{\Lambda_\ep}
|{\corr}\{ e(N_\ep),
\Lambda_\ep\cdot Z(F_\ep)\}| = 0.\vadjust{\goodbreak}
\end{equation}

To finish the proof, it suffices to prove $e(N_\ep)$ is asymptotically
uncorrelated with all linear combinations of $Z(N_\ep\cup F_\ep)$.
Specifically, defining $\Xi_\ep= (\xi_{1\ep},\ldots,\xi_{n\ep})$,
if we can show
%
%
\begin{equation}\label{corrlim}
\lim_{\ep\downarrow0}\sup_{\Lambda_\ep,\Xi_\ep}|{\corr}
\{ e(N_\ep),
\Lambda_\ep\cdot Z(F_\ep) - \Xi_\ep\cdot Z(N_\ep) \}|=0,
\end{equation}
then the theorem follows since
\[
\frac{Ee(N_\ep\cup F_\ep)^2}
{Ee(N_\ep)^2} = 1 - \sup_{\Lambda_\ep,\Xi_\ep}
\corr\{ e(N_\ep), \Lambda_\ep\cdot Z(F_\ep) - \Xi_\ep\cdot
Z(N_\ep)
\}^2.
\]
Because $e(N_\ep)$ is the error of a BLP based on $N_\ep$,
$\corr\{ e(N_\ep), \Xi_\ep\cdot Z(N_\ep)\}=0$
for all $\Xi_\ep$.
Thus, (\ref{corrlim}) follows from (\ref{corr5}) if
%
%
\begin{equation}\label{varbound}
\var\{\Lambda_\ep\cdot Z(F_\ep)\} \ll
\var\{ \Lambda_\ep\cdot Z(F_\ep) - \Xi_\ep\cdot Z(N_\ep
)\}
\end{equation}
uniformly in $\Lambda_\ep$ and $\Xi_\ep$.
There is nothing to prove if $\Xi_\ep=0$, so assume \mbox{$\Xi_\ep\ne0$}
hereafter.
Consider the Mat\'ern spectral
density $f_\alpha(\omega) = (1+\omega^2)^{-(\alpha+1)/2}$,
for which the corresponding autocovariance function is $K_\alpha
(x)\,{=}\,c_\alpha
|x|^{\alpha/2}\mathcal{K}_{\alpha/2}(|x|)$,
where $c_\alpha=
\pi^{1/2}/\{2^{\alpha/2-1}\Gamma((\alpha+1)/2)\}$ [\citet{Ste99N1},
page 31].
I will
write the subscript $\alpha$ to indicate quantities such as
variances calculated under $K_\alpha$.
Since, by~(\ref{fasym}), $f(\omega)\asymp f_\alpha(\omega)$,
(\ref{varbound}) is equivalent to
%
%
\begin{equation}
\var_\alpha\{\Lambda_\ep\cdot Z(F_\ep)\}
\ll\var_\alpha\{ \Lambda_\ep\cdot Z(F_\ep)
- \Xi_\ep\cdot Z(N_\ep)\}
\end{equation}
uniformly in $\Lambda_\ep$ and $\Xi_\ep$, which is in turn
equivalent to
%
%
\begin{equation}\label{corrbound1}
\limsup_{\ep\downarrow0}
\sup_{\Lambda_\ep,\Xi_\ep}| {\corr_\alpha}\{\Lambda_\ep
\cdot
Z(F_\ep),\Xi_\ep\cdot Z(N_\ep)\}|<1.
\end{equation}


Define $\la_{\cdot\ep}\,{=}\,\sum_{j=1}^m \la_{j\ep}$, $\tilde\la_{j\ep}
\,{=}\,\la_{j\ep}-\frac{1}{m}\la_{\cdot\ep}$, $\tilde L_\ep\,{=}\,\sum_{j=1}^m
|\tilde\la_{j\ep}|$ and $\tilde\Lambda_\ep\,{=}\,
(\tilde\la_{1\ep},\ldots,\allowbreak\tilde\la_{m\ep})$. Using\vspace*{1pt}
the series expansion for $K_\alpha$ [\citet{Ste99N1}, (15) page
32] and setting $b_\alpha=\pi/\{\Gamma(\alpha+1)\sin(\frac{1}{2}\pi
\alpha)\}$ and $S_\alpha(\ep) = -\sum_{j,k=1}^m
\tilde\la_{j\ep}\tilde\la_{k\ep} |y_j-y_k|^\alpha$,
\[
\var_\alpha\{\tilde\Lambda_\ep\cdot Z(F_\ep)\}
- b_\alpha\ep^\alpha S_\alpha(\ep) \ll\ep^2\tilde L_\ep^2.
\]
Now $S_\alpha(\ep)$ is nonnegative because $\sum_{j=1}^m \tilde
\lambda_{j\ep}
=0$, and $|x|^\alpha$ is a valid variogram for $\alpha\in(0,2)$
[\citet{Ste99N1}, page 37].
Furthermore, if $\tilde L_\ep\ne0$, $S_\alpha(\ep)/\tilde L_\ep^2$
is trivially bounded from above.
It is also uniformly bounded from below: if $S_\alpha(\ep)/\tilde
L_\ep^2$
tends to a limit along any
sequence of $\ep$ values, then there is a further subsequence along which
$\tilde\Lambda_\ep/\tilde L_\ep$
converges to some $\tilde\Lambda= (\tilde\la_{1},\ldots,\tilde\la
_{m})\ne0$
and, along this subsequence, by dominated convergence,
\[
\frac{b_\alpha S_\alpha(\ep)}{\tilde L_\ep^2}
\to\int_{-\infty}^\infty\Biggl|
\sum_{j=1}^m \tilde\lambda_j e^{i\omega x_j}\Biggr|^2 |\omega
|^{-\alpha-1}
\,d\omega> 0.
\]
Thus, no subsequence of $S_\alpha(\ep)/\tilde L_\ep^2$ can have 0 as
its limit
and
%
%
\begin{equation}\label{vartilde}
\var_\alpha\{\tilde\Lambda_\ep\cdot Z(F_\ep)\} \asymp\ep^\alpha
\tilde
L_\ep^2,
\end{equation}
which holds even if $\tilde L_\ep=0$.
Again using the series expansion for $K_\alpha$,\break
$|{\cov_\alpha}\{Z(y_0+\ep y_1),\tilde\Lambda_\ep\cdot
Z(F_\ep)\}| \ll\ep^\alpha\tilde L_\ep$,
so that
$\corr_\alpha\{Z(y_0+\ep y_1),\tilde\Lambda_\ep\cdot\break
Z(F_\ep)\} \to0$
uniformly in $\tilde\Lambda_\ep\ne0$.
Thus,
%
%
\begin{equation}\label{var1}
\var_\alpha\{\Lambda_\ep\cdot Z(F_\ep)\} \sim\la_{\cdot\ep
}^2K_\alpha(0) +
\var_\alpha\{\tilde\Lambda_\ep\cdot Z(F_\ep)\}
\end{equation}
as $\ep\downarrow0$, uniformly in $\Lambda_\ep$.
Results similar to (\ref{vartilde}) and (\ref{var1})
apply to $\Xi_\ep\cdot Z(N_\ep)$.\vspace*{1pt}

Next, define
$\xi_{\cdot\ep} = \sum_{j=1}^n \xi_{j\ep}$, $\tilde\xi_{j\ep}
=\xi_{j\ep}-\frac{1}{n}\xi_{\cdot\ep}$, $\tilde X_\ep= {\sum_{j=1}^n}
|\tilde\xi_{j\ep}|$ and $\tilde\Xi_\ep=
(\tilde\xi_{1\ep},\ldots,\tilde\xi_{n\ep})$ and
consider
%
%
\begin{eqnarray}\label{covexp}\quad
& &
\cov_\alpha\{\Lambda_\ep\cdot Z(F_\ep),\Xi_\ep\cdot Z(N_\ep)\}
\nonumber\\
& &\qquad = \la_{\cdot\ep}\xi_{\cdot\ep}K_\alpha(y_0+\ep y_1 - \ep x_1)
+ \la_{\cdot\ep}\cov_\alpha\{Z(y_0+\ep y_1),\tilde\Xi_\ep\cdot
Z(N_\ep)\}
\\
& &\qquad\quad{} + \xi_{\cdot\ep}\cov_\alpha\{Z(\ep x_1),\tilde\Lambda
_\ep\cdot
Z(F_\ep)\} + \cov_\alpha\{\tilde\Lambda_\ep\cdot Z(F_\ep),\tilde
\Xi_\ep\cdot
Z(N_\ep)\}.
\nonumber
\end{eqnarray}
Since $K_\alpha$ has a bounded second derivative outside of a
neighborhood of
the origin, it is straightforward to obtain the following bounds:\vspace*{2pt}
\begin{eqnarray*}
|\la_{\cdot\ep}\xi_{\cdot\ep}K_\alpha(y_0+\ep y_1 - \ep x_1) -
\la_{\cdot\ep}\xi_{\cdot\ep}K_\alpha(y_0)| &\ll&\ep|\la_{\cdot
\ep}||
\xi_{\cdot\ep}|,
\\
|\la_{\cdot\ep}\cov_\alpha\{Z(y_0+\ep y_1),\tilde\Xi_\ep\cdot
Z(N_\ep)\}|
&\ll&\ep|\la_{\cdot\ep}|\tilde X_\ep,
\\
|\xi_{\cdot\ep}\cov_\alpha\{Z(\ep x_1),\tilde\Lambda_\ep\cdot
Z(F_\ep)\}| &\ll&\ep|\xi_{\cdot\ep}|\tilde L_\ep\vspace*{2pt}
\end{eqnarray*}
and
\[
|{\cov_\alpha}\{\tilde\Lambda_\ep\cdot Z(F_\ep),\tilde\Xi_\ep
\cdot
Z(N_\ep)\}| \ll\ep^2\tilde L_\ep\tilde X_\ep
\]
as $\ep\downarrow0$.
Applying these bounds to (\ref{covexp}) gives
%
%
\begin{eqnarray}\label{b5}
& &
|{\cov_\alpha}\{\Lambda_\ep\cdot Z(F_\ep),\Xi_\ep\cdot Z(N_\ep)\} -
\la_{\cdot\ep}\xi_{\cdot\ep}K_\alpha(y_0)|
\nonumber\\[-8pt]\\[-8pt]
& &\qquad \ll\ep|\la_{\cdot\ep}||
\xi_{\cdot\ep}| + \ep|\la_{\cdot\ep}|\tilde X_\ep+ \ep|\xi
_{\cdot\ep}|\tilde
L_\ep+ \ep^2\tilde L_\ep\tilde X_\ep.
\nonumber
\end{eqnarray}
Now, from (\ref{var1}),
%
%
\begin{eqnarray}\label{Kalpha}\quad
& & \frac{|\la_{\cdot\ep}\xi_{\cdot\ep}K_\alpha(y_0)|}{\sqrt
{\var_\alpha\{
\Lambda_\ep\cdot Z(F_\ep)\}\var_\alpha\{\Xi_\ep\cdot Z(N_\ep)\}}}
\nonumber\\
& &\qquad \sim
\frac{|\la_{\cdot\ep}\xi_{\cdot\ep}K_\alpha(y_0)|}{\sqrt{\la
_{\cdot\ep}^2
K_\alpha(0)+\var_\alpha\{\tilde\Lambda_\ep\cdot Z(F_\ep)\}}
\sqrt{\xi_{\cdot\ep}^2 K_\alpha(0)+\var_\alpha\{\tilde\Xi_\ep
\cdot
Z(N_\ep)\}}}
\\
& &\qquad \le\frac{K_\alpha(y_0)}{K_\alpha(0)},
\nonumber
\end{eqnarray}
which is in $(0,1)$ for all $y_0\ne0$.
And, since $\alpha<2$,
\begin{eqnarray*}
\frac{\ep|\la_{\cdot\ep}||
\xi_{\cdot\ep}| + \ep|\la_{\cdot\ep}|\tilde X_\ep+ \ep|\xi
_{\cdot\ep}|\tilde
L_\ep+ \ep^2\tilde L_\ep\tilde X_\ep}{\sqrt{\la_{\cdot\ep
}^2+\ep^\alpha
\tilde L_\ep^2}\sqrt{\xi_{\cdot\ep}^2+\ep^\alpha
\tilde X_\ep^2}} \to0,\vadjust{\goodbreak}
\end{eqnarray*}
which, together with (\ref{var1}), (\ref{b5}) and (\ref{Kalpha}),
proves (\ref{corrbound1}) and hence (\ref{corrlim}) and the theorem.

To prove that Theorem \ref{thm1} also applies to ordinary kriging, note
that by setting $\beta=\frac{1}{2}$, (\ref{phibound}) and (\ref{eta-tilde})
together with (\ref{rate}) imply $\|\tilde\eta_\ep\|_f^2\sim\|\eta
_\ep\|_f^2$
as $\ep\downarrow0$.
Since $\tilde\eta_\ep$ corresponds to the error of a linear unbiased
predictor under the constant mean model, we have that the mean squared
error of the ordinary kriging predictor based on $Z(N_\ep)$ is at
least $\|\eta_\ep\|_f^2$ and at most $\|\tilde\eta_\ep\|_f^2$, so that
if (\ref{main}) holds for the simple kriging predictor it also holds
for the
ordinary kriging predictor.

\subsection{\texorpdfstring{Proof of Theorem \protect\ref{thm2}}{Proof of Theorem 2}}\label{sec53}

Restricting $N_\ep$ to one point and $F_\ep$ to 2 allows us to make
use of
Lemma \ref{lem1} to prove (\ref{main}).
Setting $y_2=0$ simplifies the calculations without changing any essential
details.
Specifically, defining $V(x) = K(0)-K(x)$, we will show that
%
%
\begin{eqnarray}\label{t2vec}\quad
W'_\ep& = & (W_{\ep1},W_{\ep2},W_{\ep3},W_{\ep4}) \nonumber\\[-8pt]\\[-8pt]
& = &
\biggl( \frac{Z(0)-Z(\ep x_1)}{\sqrt{V(\ep x_1)}},Z(\ep
x_1),Z(y_0+\ep y_1),
\frac{Z(y_0)-Z(y_0+\ep y_1)}{\sqrt{V(\ep y_1)}}\biggr)'
\nonumber\hspace*{-25pt}
\end{eqnarray}
has limiting covariance matrix of the form given in (\ref{lemma-limit}),
from which Theorem~\ref{thm2} readily follows.

Let us consider the easier parts of the proof first.
Independent of $\ep$,
the variances of the elements of $W_\ep$ are $1,K(0),K(0)$ and 1, respectively.
Since~$Z$ has a spectral density, $K$ is continuous and $|K(y)| < K(0)$
for all $y\ne0$.
Thus, $\cov(W_{\ep2},W_{\ep3}) \to K(y_0)$ as $\ep\downarrow0$, and
the $2\times2$ matrix with $K(0)$ on the diagonals and $K(y_0)$
elsewhere is
positive definite.
Thus, it suffices to show that the other offdiagonal elements of the
covariance matrix of $W_\ep$ tend to $0$ as $\ep\downarrow0$.
First, $\cov(W_{\ep1},W_{\ep2}) = \frac{1}{2}\sqrt{V(\ep x_1)/K(0)}
\to0$ as $\ep\downarrow0$.
Similarly, $\cov(W_{\ep3},W_{\ep4})\to0$ as $\ep\downarrow0$.

Now consider $\cov(W_{\ep1},W_{\ep3})$.
We have
\[
\cov\{Z(0)-Z(\ep x_1),Z(y_0+\ep y_1)\} =
\int_{\R^2} e^{-i\omega\cdot(y_0+\ep y_1)}(1-e^{i\ep\omega\cdot
x_1})f(\omega)
\,d\omega,
\]
so that for $D(T) = \{\omega\dvtx |\omega\cdot x_1| \le T \}$,
%
%
\begin{eqnarray}\label{corr13}
\cov(W_{\ep1},W_{\ep3})^2
&\le&\frac{\{\int_{\R^2} |1-e^{i\ep\omega\cdot x_1}
| f(\omega)\,d\omega\}^2}{K(0)\int_{\R^2} |1-e^{i\ep
\omega\cdot x_1}
|^2 f(\omega)\,d\omega}\nonumber\\[-0.5pt]
&\le&\frac{2\{\int_{D(T)} |1-e^{i\ep\omega\cdot x_1}
| f(\omega)\,d\omega\}^2}{K(0)\int_{D(T)}
|1-e^{i\ep\omega\cdot x_1}
|^2 f(\omega)\,d\omega} \\[-0.5pt]
&&{}+ \frac{2\{\int_{D(T)^c}
|1-e^{i\ep\omega\cdot
x_1} | f(\omega)\,d\omega\}^2}{K(0)\int_{D(T)^c}
|1-e^{i\ep\omega\cdot x_1} |^2 f(\omega)\,d\omega}
\nonumber
\end{eqnarray}
for all $T$ sufficiently large (to guarantee $\int_{D(T)}
|1-e^{i\ep\omega\cdot x_1}
|^2 f(\omega)\,d\omega> 0$).
Because $\ep^{-1}|1-e^{i\ep\omega\cdot x_1} | \le|\omega
\cdot x_1|$
and $\ep^{-1}|1-e^{i\ep\omega\cdot x_1} | \to|\omega
\cdot x_1|$
as $\ep\downarrow0$, by dominated convergence,
%
%
\begin{equation}
\lim_{\ep\downarrow0}\label{eplim1}
\frac{\{\int_{D(T)} |1-e^{i\ep\omega\cdot x_1}
| f(\omega)\,d\omega\}^2}{\int_{D(T)}
|1-e^{i\ep\omega\cdot x_1} |^2 f(\omega)\,d\omega}
= \frac{\{\int_{D(T)}|\omega\cdot x_1|f(\omega)\,d\omega\}^2}
{\int_{D(T)} |\omega\cdot x_1|^2 f(\omega)\,d\omega}.
\end{equation}
By the Cauchy--Schwarz inequality,
%
%
\begin{eqnarray}\label{eplim2}
&&\frac{\{\int_{D(T)^c}
|1-e^{i\ep\omega\cdot
x_1} | f(\omega) \,d\omega\}^2}{\int_{D(T)^c}
|1-e^{i\ep\omega\cdot x_1} |^2 f(\omega) \,d\omega}\nonumber\\[-0.5pt]
&&\qquad \le \frac{\int_{D(T)^c}
|1-e^{i\ep\omega\cdot x_1} |^2 f(\omega) \,d\omega
\int_{D(T)^c} f(\omega) \,d\omega}{\int_{D(T)^c}
|1-e^{i\ep\omega\cdot x_1}|^2 f(\omega) \,d\omega}
\\[-0.5pt]
&&\qquad = \int_{D(T)^c} f(\omega)\,d\omega.
\nonumber
\end{eqnarray}
From (\ref{corr13})--(\ref{eplim2}), we will have $\cov(W_{\ep
1},W_{\ep
3})\to0$ as $\ep\downarrow0$ if the right-hand sides of (\ref
{eplim1}) and
(\ref{eplim2}) tend to 0 as $T\to\infty$.
The integrability of $f$ implies $\int_{D(T)^c} f(\omega)\,d\omega\to
0$ as
$T\to\infty$, so consider\vspace*{1pt} the right-hand side of (\ref{eplim1}).
Let~$A$ be the $2\times2$ matrix with first row given by $x_1$,
orthogonal\vspace*{1pt} rows and determinant of 1 and set $v = (v_1, v_2)' = A\omega$.
Define $\bar f(v_1) = \int_{-\infty}^{\infty} f( A^{-1} v
) \,dv_2$.
Up to a linear rescaling, $\bar f$ is
the spectral density of the process $Z$ along the $x_1$ direction,
so it is integrable.
In addition, because $Z$ is not mean square differentiable in
any direction, $\int_0^\infty v_1^2\bar f(v_1) \,dv_1 = \infty$.
Then for any even function~$g$, $\int_{D(T)}g(\omega\cdot
x_1)f(\omega)
\,d\omega= 2\int_0^T g(v_1)\bar f(v_1) \,dv_1$, so that for $0<S<T$,
%
%
\begin{eqnarray}\label{0ST}
\frac{\{\int_{D(T)}|\omega\cdot x_1|f(\omega)\,d\omega\}^2}
{\int_{D(T)} |\omega\cdot x_1|^2 f(\omega)\,d\omega}
& = & \frac{\{ \int_0^T v_1\bar f(v_1) \,dv\}^2}
{ \int_0^T v_1^2\bar f(v_1) \,dv_1} \nonumber\\[-8pt]\\[-8pt]
& = & \frac{\{ \int_0^S v_1\bar f(v_1) \,dv_1 +
\int_S^T v_1\bar f(v_1) \,dv_1\}^2}
{ \int_0^T v_1^2\bar f(v_1) \,dv_1}.
\nonumber
\end{eqnarray}
If we can show that
%
%
\begin{equation}\label{lim0ST}
\lim_{S\to\infty}\lim_{T\to\infty} \frac{\{ \int_0^S
v_1\bar f(v_1)
\,dv_1
+ \int_S^T v_1\bar f(v_1) \,dv_1\}^2} { \int_0^T v_1^2\bar
f(v_1)
\,dv_1} = 0,
\end{equation}
then the right-hand side of (\ref{eplim1}) will tend to 0 as $T\to
\infty$.
To prove (\ref{lim0ST}), expand the square in the numerator and consider
each term separately.
First, by the Cauchy--Schwarz inequality,
\[
\lim_{T\to\infty}
\frac{\{ \int_0^S v_1\bar f(v_1) \,dv_1\}^2}
{ \int_0^T v_1^2\bar f(v_1) \,dv_1}
= \lim_{T\to\infty}\frac{ \int_0^S v_1^2\bar f(v_1) \,dv_1 \int_0^S
\bar f(v_1) \,dv_1}
{ \int_0^T v_1^2\bar f(v_1) \,dv_1} = 0.
\]
Again by the Cauchy--Schwarz inequality,
\[
\frac{\{
\int_S^T v_1\bar f(v_1) \,dv_1\}^2}
{ \int_0^T v_1^2\bar f(v_1) \,dv_1}
\le\frac{ \int_S^T v_1^2\bar f(v_1) \,dv_1 \int_S^T \bar f(v_1) \,dv_1}
{ \int_0^T v_1^2\bar f(v_1) \,dv_1} \\
\le\int_S^T \bar f(v_1) \,dv_1,
\]
which tends to 0 when one takes $\lim_{S\to\infty}\lim_{T\to\infty
}$ since
$\bar f$ is integrable.
Finally,
\begin{eqnarray*}
\frac{ \int_0^S v_1\bar f(v_1) \,dv_1 \int_S^T v_1\bar f(v_1) \,dv_1}
{\int_0^T v_1^2\bar f(v_1) \,dv_1}
&\le& \frac{ \int_0^S v_1\bar f(v_1) \,dv_1 \int_S^T v_1\bar
f(v_1)
\,dv_1} { \int_S^T v_1^2\bar f(v_1) \,dv_1} \\
&\le& \frac{ \int_0^S v_1\bar f(v_1) \,dv_1}{S}\cdot
\frac{\int_S^T v_1\bar f(v_1) \,dv_1} {\int_S^T v_1\bar f(v_1)
\,dv_1} \\
&\le& \frac{1}{S}\int_0^{S^{1/2}} v_1\bar f(v_1) \,dv_1 +
\frac{1}{S}\int_{S^{1/2}}^S v_1\bar f(v_1) \,dv_1 \\
&\le& \frac{1}{S^{1/2}} \int_0^{S^{1/2}}\bar f(v_1) \,dv_1
+ \int_{S^{1/2}}^S \bar f(v_1) \,dv_1,
\end{eqnarray*}
which tends to 0 as $S\to\infty$, and (\ref{lim0ST}) follows.
Thus, $\cov(W_{\ep1},W_{\ep3})\to0$ as $\ep\downarrow0$.
Similarly, $\cov(W_{\ep2},W_{\ep4})\to0$ as $\ep\downarrow0$.

We will need the following lemma to handle $\cov(W_{\ep1},W_{\ep4})$:
\begin{lemma}\label{lem2}
If $Z$ is not mean square differentiable in the direction $x$, then
\[
\lim_{\ep\downarrow0}\frac{\ep^2}{V(\ep x)} = 0.
\]
\end{lemma}

To prove the lemma, first note that the assumption on $Z$ is equivalent
to
%
%
\begin{equation}\label{no2}
\int_{\R^2} |\omega\cdot x|^2 f(\omega) \,d\omega= \infty.
\end{equation}
If $\limsup_{\ep\downarrow0}\frac{\ep^2}{V(\ep x)} > 0$, then there
must exist some sequence $\ep_n\downarrow0$ along which
$\lim_{n\to\infty}\frac{V(\ep_n x)}{\ep_n^2} = C$ for some finite
$C$, or
\[
\lim_{n\to\infty}\int_{\R^2} \frac{|1-e^{i\ep_n\omega\cdot
x}|^2}{\ep_n^2}
f(\omega) \,d\omega= C.
\]
But for any finite $T$, by dominated convergence,
\begin{eqnarray*}
C & = & \lim_{n\to\infty}\int_{\R^2} \frac{|1-e^{i\ep_n\omega
\cdot x}|^2}
{\ep_n^2} f(\omega) \,d\omega\\
& \ge& \lim_{n\to\infty}\int_{|\omega|<T} \frac{|1-e^{i\ep
_n\omega\cdot
x}|^2}{\ep_n^2} f(\omega) \,d\omega\\
& = & \int_{|\omega|<T} |\omega\cdot x|^2 f(\omega) \,d\omega
\end{eqnarray*}
for all $T$, which contradicts (\ref{no2}), and the lemma is proven.

Consider
\begin{eqnarray*}
& & \cov\{Z(0)-Z(\ep x_1),Z(y_0)-Z(y_0+\ep y_1)\} \\
& &\qquad = \int_{\R^2} e^{-i\omega\cdot y_0}(1-e^{i\ep\omega\cdot
x_1})
(1-e^{-i\ep\omega\cdot y_1}) f(\omega) \,d\omega.
\end{eqnarray*}
Define $f_1(\omega) = \min(f(\omega),1)$, and write $\cov_1$ to indicate
covariances calculated under the spectral density $f_1$.
Then (\ref{f-cond}) and $f$ integrable
imply that $f(\omega) = f_1(\omega)$ outside some bounded set, and it easily
follows that
%
%
\begin{eqnarray}\label{f1}
& & \cov\{Z(0)-Z(\ep x_1),Z(y_0)-Z(y_0+\ep y_1)\} \nonumber\\[-8pt]\\[-8pt]
& &\qquad = \cov_1\{Z(0)-Z(\ep x_1),Z(y_0)-Z(y_0+\ep y_1)\} + O(\ep^2).
\nonumber
\end{eqnarray}
Because $Z$ is not mean square differentiable in any direction, Lemma
\ref{lem2}
implies the $O(\ep^2)$ remainder in (\ref{f1}) makes no
contribution to $\lim_{\ep\downarrow0}\cov(W_{\ep1},W_{\ep4})$.

We proceed by rotating coordinates so that one of the frequency axes points
in the direction of $y_0$.
Specifically,
let $B$ be the $2\times2$ orthogonal matrix with determinant 1 and
first row equal to $y_0$ and set $\tau= (\tau_1, \tau_2)' = B\omega$.
Then, defining $H_\ep(\tau_1,\tau_2) =
(1-e^{i\ep(B^{-1}\tau)\cdot x_1})
(1-e^{-i\ep(B^{-1}\tau)\cdot y_1})$,
\begin{eqnarray*}
& & \cov_1\{Z(0)-Z(\ep x_1),Z(y_0)-Z(y_0+\ep y_1)\} \\
& &\qquad = \int_{\R^2} e^{-i\tau_1}H_\ep(\tau_1,\tau_2)
f_1(B^{-1}\tau) \,d\tau\\
& &\qquad = \int_\R\sum_{k=-\infty}^{\infty} \int_{2\pi k}^{2\pi(k+1)}
e^{-i\tau_1}H_\ep(\tau_1,\tau_2)
f_1(B^{-1}\tau) \,d\tau_1
\,d\tau_2.
\end{eqnarray*}
Define the function $g$
on $\R^2$ by, for $2\pi k\le\tau_1 < 2\pi(k+1)$, $g(B^{-1}\tau) =
1 -
f_1(B^{-1}(2\pi k, \tau_2)')/f_1(B^{-1}\tau)$ if
$f_1(B^{-1}\tau) >0$ and 0 otherwise.
We have
%
%
\begin{eqnarray}\label{splitk}
& & \cov_1\{Z(0)-Z(\ep x_1),Z(y_0)-Z(y_0+\ep y_1)\} \nonumber\\
& &\qquad = \int_\R\sum_{k=-\infty}^{\infty} f_1\biggl( B^{-1}\pmatrix
{2\pi k\cr
\tau_2}\biggr) \int_{2\pi k}^{2\pi(k+1)}
e^{-i\tau_1}H_\ep(\tau_1,\tau_2) \,d\tau_1
\,d\tau_2\\
& &\qquad\quad{} +
\int_{\R^2} e^{-i\omega\cdot y_0}(1-e^{i\ep\omega\cdot
x_1})
(1-e^{-i\ep\omega\cdot y_1}) f_1(\omega)g(\omega)
\,d\omega.
\nonumber
\end{eqnarray}
By (\ref{f-cond}), $g(\omega)\to0$ as $\omega\to\infty$.
Thus, given $\delta>0$, we can find $T<\infty$ such that $g(\omega
)<\delta$
for $|\omega|>T$.
Then
\begin{eqnarray*}
& &
\biggl| \int_{\R^2} e^{-i\omega\cdot y_0}(1-e^{i\ep\omega
\cdot x_1})
(1-e^{-i\ep\omega\cdot y_1}) f_1(\omega)g(\omega)
\,d\omega\biggr|
\\
& &\qquad \le\ep^2 \int_{|\omega|\le T} |\omega\cdot x_1||\omega\cdot y_1|
f_1(\omega)|g(\omega)| \,d\omega\\
& &\qquad\quad{} + 4\delta\int_{|\omega|>T}
|1-e^{i\ep\omega\cdot x_1}|
|1-e^{-i\ep\omega\cdot y_1}| f_1(\omega) \,d\omega.
\end{eqnarray*}
By the Cauchy--Schwarz inequality and $f_1\le f$, $\int_{|\omega|>T}
|1-e^{i\ep\omega\cdot x_1}|
|1-\break e^{-i\ep\omega\cdot y_1}|\times f_1(\omega) \,d\omega\le
\sqrt{V(\ep x_1)V(\ep y_1)}$,
which, together with Lemma \ref{lem2}, implies
%
%
\begin{equation}\label{4delta}
\limsup_{\ep\downarrow0}
\frac{|
\int_{\R^2} e^{-i\omega\cdot y_0}(1-e^{i\ep\omega\cdot
x_1})
(1-e^{-i\ep\omega\cdot y_1}) f_1(\omega)g(\omega)
\,d\omega|}
{\sqrt{V(\ep x_1)V(\ep y_1)}} \le4\delta.
\end{equation}
Since $\delta$ is arbitrary, this $\limsup$ must in fact be 0.

Now return to the the first term on the right-hand side of (\ref{splitk}).
Integrating by parts,
\begin{eqnarray*}
\int_{2\pi k}^{2\pi(k+1)} e^{-i\tau_1}H_\ep(\tau_1,\tau_2)
\,d\tau_1
&=& iH_\ep\bigl(2\pi(k+1),\tau_2\bigr)- iH_\ep(2\pi k,\tau_2)\\
&&{} -
i\int_{2\pi k}^{2\pi(k+1)} e^{-i\tau_1}\,\frac{\partial}{\partial
\tau_1}
H_\ep(\tau_1,\tau_2) \,d\tau_1.
\end{eqnarray*}
There exists finite $C$ independent of $\ep$ and $\tau$ such that
\[
\biggl| \frac{\partial}{\partial\tau_1}
H_\ep(\tau_1,\tau_2)\biggr| \le C\ep\bigl\{\bigl|1-e^{i\ep
(B^{-1}\tau)\cdot
x_1}\bigr|+\bigl|1-e^{-i\ep(B^{-1}\tau)\cdot y_1}\bigr|\bigr\},
\]
which implies
%
%
\begin{eqnarray}\label{Hbound}
& & \biggl|\int_{2\pi k}^{2\pi(k+1)} e^{-i\tau_1}H_\ep(\tau_1,\tau
_2) \,d\tau_1\biggr|
\nonumber\\[-8pt]\\[-8pt]
& &\qquad \le4\pi C\ep
\bigl\{\bigl|1-e^{i\ep(B^{-1}\tau)\cdot
x_1}\bigr|+\bigl|1-e^{-i\ep(B^{-1}\tau)\cdot y_1}\bigr|\bigr\}.
\nonumber
\end{eqnarray}
We can choose $T$ finite so that if
$2\pi k\le\tau_1 \le2\pi(k+1)$, then $f_1( B^{-1}(2\pi k,\allowbreak \tau
_2)')
\le2 f( B^{-1}\tau)$ whenever $|\tau| > T$.
Applying this result and (\ref{Hbound})
to the first term on the right-hand side of
(\ref{splitk}) and changing variables back to $\omega= B^{-1}\tau$,
we get
%
%
\begin{eqnarray}\label{laststep}
& & \Biggl| \int_\R\sum_{k=-\infty}^{\infty} f_1\biggl( B^{-1}
\pmatrix{2\pi k\cr\tau_2}\biggr) \int_{2\pi k}^{2\pi(k+1)}
e^{-i\tau_1}H_\ep(\tau_1,\tau_2)
\,d\tau_1 \,d\tau_2\Biggr| \nonumber\\
& &\qquad \le8\pi C\ep\int_{\R^2} f(\omega)
\{|1-e^{i\ep\omega\cdot
x_1}|+|1-e^{-i\ep\omega\cdot y_1}|\} \,d\omega
+O(\ep^2)
\\
& &\qquad \le8\pi C\ep\bigl\{\sqrt{V(\ep x_1)}+\sqrt{V(\ep y_1)}\bigr\}
\sqrt{\int_{\R^2} f(\omega) \,d\omega} + O(\ep^2),
\nonumber
\end{eqnarray}
where the last step uses the Cauchy--Schwarz inequality.
From Lemma \ref{lem2} and~(\ref{laststep}), it follows that
\[
\limsup_{\ep\downarrow0}
\frac{| {\int_\R}
\sum_{k=-\infty}^{\infty} f_1( B^{-1}{2\pi
k\choose\tau_2}) \int_{2\pi k}^{2\pi(k+1)}
e^{-i\tau_1}H_\ep(\tau_1,\tau_2)
\,d\tau_1 \,d\tau_2|}
{\sqrt{V(\ep x_1)V(\ep y_1)}} = 0.
\]
Together with (\ref{splitk}) and (\ref{4delta}), this limit implies
\[
\limsup_{\ep\downarrow0}
\frac{\cov_1 \{Z(0)-Z(\ep x_1),Z(y_0)-Z(y_0+\ep y_1)\}}
{\sqrt{V(\ep x_1)V(\ep y_1)}} = 0,
\]
which together with (\ref{f1}) and Lemma \ref{lem2}, implies
$\lim_{\ep\downarrow0} \cov\{W_{\ep1},W_{\ep4}\} = 0$.

Theorem \ref{thm2} applies to ordinary kriging as well.
Specifically, $Z(\ep x_1)$ is an asymptotically
optimal linear predictor of $Z(0)$ based on $Z(N_\ep\cup F_\ep)$
when the mean of $Z$ is assumed to be 0,
so since it
is a linear unbiased predictor when the mean is an unknown
constant, $Z(\ep x_1)$
must also be asymptotically optimal with respect to this
more restricted class of predictors.


\section*{Acknowledgment}

The author thanks Steven Lalley for help with the proof of Theorem
\ref{thm2}.


%

\printaddresses

\end{document}